\documentclass[a4paper,11pt,draft]{article}
\usepackage[english]{babel}
\usepackage{theorem}
\usepackage{amssymb}
\usepackage{amsfonts}
\usepackage{amsmath}
\usepackage[all]{xy}
\usepackage{geometry}
\usepackage{ae,aecompl}
\usepackage{eurosym}
\usepackage{verbatim}
\usepackage{mathrsfs}

\newtheorem{thm}{Theorem}[section]
\newtheorem{prop}[thm]{Proposition}
\newtheorem{lem}[thm]{Lemma}
\newtheorem{cor}[thm]{Corollary}

\newtheorem{remark}[thm]{Remark}

\makeatletter
\begingroup
\gdef\th@upshape{\normalfont
  \def\@begintheorem##1##2{%
        \item[\hskip\labelsep \theorem@headerfont ##1\ ##2.]}%
\def\@opargbegintheorem##1##2##3{%
   \item[\hskip\labelsep \theorem@headerfont ##1\ ##2\ (##3).]}}
\endgroup
\makeatother

\theoremstyle{upshape}

\newtheorem{eg}[thm]{Example}
\newtheorem{rem}[thm]{Remark}
\newtheorem{defn}[thm]{Definition}
\newtheorem{conj}[thm]{Conjecture}

\date{}

\def\bal#1\eal{\begin{align}#1\end{align}}
\def\bas#1\eas{\begin{align*}#1\end{align*}}
\def\bit{\begin{itemize}}
\def\eit{\end{itemize}}
\def\bet{\begin{enumerate}}
\def\eet{\end{enumerate}}
\def\ed{\end{document}}

\def\a{\alpha}
\def\b{\beta}
\def\g{\gamma}
\def\d{\delta}
\def\e{\varepsilon}
\def\f{\varphi}
\def\k{\kappa}
\def\l{\lambda}
\def\Lam{\Lambda}

\def\r{\rho}
\def\s{\sigma}

\def\t{\tau}
\def\w{\omega}

\def\Om{\Omega}

\def\del{\partial}
\def\adel{\ol{\partial}}
\def\DEL{\Delta}
\def\G{\Gamma}

\def\bC{{\mathbb C}}
\def\bN{{\mathbb N}}

\def\bR{{\mathbb R}}
\def\bZ{{\mathbb Z}}

\def\bT{\mathbb{T}}

\def\A{{\cal A}}
\def\B{{\cal B}}
\def\C{{\cal C}}
\def\D{{\cal D}}
\def\E{{\cal E}}
\def\F{{\cal F}}
\def\H{{\cal H}}

\def\exd{\mathrm{d}}
\def\demo{\noindent \emph{\textbf{Proof}.\ }~}
\def\dim{\mathrm{dim}}

\def\haar{\mathrm{\bf h}}
\def\unit{\mathrm{U}}
\def\counit{\mathrm{C}}

\def\id{\mathrm{id}}
\def\ker{\mathrm{ker}}
\def\proj{\mathrm{proj}}

\def\spn{\mathrm{span}}
\def\th{^\mathrm{th}}

\def\vol{\mathrm{vol}}

\def\hol{^{(1,0)}}
\def\ahol{^{(0,1)}}

\def\inv{^{-1}}

\def\by{\times}
\def\coby{\, \square_{H}}
\def\oby{\otimes}

\def\wed{\wedge}
\def\sseq{\subseteq}

\def\wt{\widetilde}
\def\wh{\widehat}
\def\ol{\overline}

\def\la{\left\langle}
\def\\la{\left\langle}
\def\ra{\right\rangle}
\def\>{\right\rangle}
\def\bs{\backslash}
\def\mto{\mapsto}

\def\qed{\hfill\ensuremath{\square}\par}
\def\cf3{\bC_q[F_3]}

\def\csu2{\bC_q[SU_2]}

\def\cs1{\bC P^{1}}
\def\cccp1{\bC_q[\bC P^{1}]}
\def\ccp2{\bC P^{2}}
\def\cp2{\bC_q[\bC P^{2}]}
\def\cpn{\bC_q[\bC P^{n}]}

\def\ccpn{\bC P^{n}}

\def\qf3{\bC_q[F_3]}
\def\wqf3{\Om^1_q[F_3]}
\def\qsu3{\bC_q[SU_3]}
\def\wqsu3{\Om^1_q[SU_3]}

\def\usl2{\mathcal{U}(\mathfrak{sl}(2))}

\def\ws2{\Om^1_q(S^2)}
\def\wsu2{\Om^1_q(SU_2)}

\def\hol{^{(1,0)}}
\def\ahol{^{(0,1)}}

\def\n2{_{(-2)}}
\def\m2{_{(-2)}}
\def\m1{_{(-1)}}
\def\0{_{(0)}}
\def\1{_{(1)}}
\def\2{_{(2)}}
\def\3{_{(3)}}
\def\4{_{(4)}}
\def\5{_{(5)}}
\def\hol{^{(1,0)}}
\def\ahol{^{(0,1)}}

\def\mm{{}_M \hspace{-.030cm}\mathrm{Mod}}

\def\gmmm{{}^{G}_M \hspace{-.030cm}\mathrm{Mod}_M}

\def\sgmm{{}^{G}_M \hspace{-.030cm}{\mathrm{mod}}_0}
\def\smh{{^H\!\mathrm{mod}}_0}
\def\lgmmm{{}^{G}_M \hspace{-.030cm}{\mathrm{Mod}}_0}
\def\lmhm{{^H\!\mathrm{Mod}}_0}

\def\alg{algebra~}
\def\algn{algebra}
\def\algs{algebras~}

\def\hk{Heckenberger--Kolb~}

\def\ncg{noncommutative geometry~}

\def\nc{noncommutative~}

\def\onb{orthonormal basis~}

\def\qhs{quantum homogeneous space~}
\def\qhsn{quantum homogeneous space}

\def\iff{if and only if~}

\def\st{such that~}

\def\wrt{with respect to~}

\DeclareMathOperator{\dt}{det}

\def\mpr{maximal prolongation~}

\def\uqsl2{U_q(\frak{sl_2)}}
\def\tuqsl2{\wt{U}_q(\frak{sl}_2)}
\def\tu1sl2{\wt{U}_1(\frak{sl}_2)}

\def\sgmm{{}^{G}_M \hspace{-.030cm}{\mathrm{mod}}_0}

\title{{\bf Noncommutative K\"ahler  Structures on Quantum Homogeneous Spaces}}
\author{R\'{e}amonn \'{O} Buachalla\footnote{The paper was supported by funds allocated to the implementation of the international co-funded project in the years 2014-2018 3038/7.PR/2014/2 and by the grants FP7-PEOPLE-2012-COFUND-600415 and GACR P201/12/G028.}}

\begin{document}

\maketitle

\begin{abstract}
Building on the theory of noncommutative complex structures, the notion  of a  noncommutative K\"ahler structure is introduced. In the quantum homogeneous space case many of the fundamental results of classical K\"ahler geometry are shown to follow from the existence of such a structure, allowing for the definition of noncommutative Lefschetz, Hodge, K\"ahler--Dirac, and Laplace operators. Quantum projective space, endowed with its \hk calculus, is taken as the motivating example. The general theory is then used to show that the calculus has cohomology groups of at least classical dimension.
\end{abstract}

\section{Introduction}

One of the most exciting new trends in \ncg is the search for a theory of noncommutative complex geometry \cite{KLVSCP1, BS, MMF2}. It is motivated by the appearance of noncommutative complex structures in a number of areas of \nc geometry, such as the construction of spectral triples for quantum groups  \cite{Krah, DDCPN, EBSMSTrip},  geometric representation theory for quantum groups \cite{JS, Maj, KLVSCP1, KKCPN},  the interaction of noncommutative geometry and  noncommutative projective  algebraic geometry \cite{KLVSCP1, KKCPN, BS}, the Baum--Connes conjecture for quantum groups \cite{V,VY}, and the application of topological algebras to quantum group theory \cite{AP1,AP2}.  While there have been a number of occurrences of K\"ahler phenomena in the literature, the question of whether metrics have a role to play in noncommutative complex geometry remains largely unexplored.  Given the richness and beauty of classical K\"ahler geometry, the idea that some of its structure might generalise to the \nc setting is an enticing one.

To date there has been just one proposed framework for  \nc K\"ahler geometry \cite{FGR}.  It uses a Riemannian, as opposed to spin, approach to spectral triples, and takes as its motivating example the \nc torus. In this paper we instead take the quantum flag manifolds as a motivating family of examples, and adopt an approach based on Woronowicz's notion of a differential calculus \cite{Wor}.  Classically the flag manifolds play a central role in  K\"ahler and parabolic geometry \cite{ParbolicG}, and as such, their $q$-deformations serve as outstanding candidates for noncommutative K\"ahler spaces.  Moreover, the subfamily of irreducible quantum flag manifolds comes endowed with a  differential calculus, the \hk calculus, which is uniquely characterised by a simple set of natural axioms \cite{HK, HKdR}. These calculi have already served as the motivating examples for the theory of  \nc complex structures \cite{BS, KLVSCP1, MMF2}.

Metric phenomena have appeared a number of times in the literature on the \ncg of the quantum flag manifolds. A Hodge map for the Podle\'s sphere was defined by Majid in \cite{Maj}, and the induced Laplace and Dirac operators studied. Hodge maps on the Podle\'s sphere and $\bC_q[\bC P^2]$ were examined by Landi, Zampini, and D'Andrea in the series of papers  \cite{DLASD,LZ,Z}. A K\"ahler--Dirac operator for the irreducible quantum flag manifolds was introduced by Kr\"ahmer in \cite{Krah} and used to give a commutator presentation of the \hk calculus. This operator was reconstructed by D\c{a}browski, Landi, and D'Andrea in \cite{DDLCP2,DDCPN} for the special case of $\cpn$ and was shown to satisfy the properties of a spectral triple. Finally, a direct $q$-deformation of the  K\"ahler form of $\bC P^2$ was constructed in \cite{DLASD}.

Inspired by such phenomena, this paper introduces a general framework for \nc K\"ahler geometry on quantum homogeneous spaces and applies it to  quantum projective space.  The manner in which this is done has three main sources of inspiration: The first is Majid's frame bundle approach to noncommutative geometry  \cite{qqguage,Maj1,Maj}  which also underpins the author's earlier papers \cite{MMF1,MMF2}. The second is Kustermanns, Murphy, and Tuset's approach to noncommutative Hodge theory \cite{KMT}, and the third is the  presentation of classical K\"ahler geometry found in \cite{Weil} and \cite{HUY}, which is both global and algebraic in style.

The major obstacle to formulating a coherent construction of metrics in the noncommutative setting is that the classical extension of metrics from $1$-forms to higher forms does not easily generalise. The classical extension uses anti-commutativity of forms in a fundamental way, while the multiplicative relations for a differential calculus over a \nc \alg will in general be much more badly behaved.  In certain cases, such as bicovariant calculi \cite{Wor}, or braided complex structures \cite{MMF5}, one can formulate a braided generalisation of the classical construction \cite{K}. However, in practice  the metrics produced are not ideal \cite{MMF5}.

For  K\"ahler manifolds, however,  we show that it is possible to reverse the usual order of construction and build a Hodge map from a K\"ahler form and then use this to extend the metric. Adopting this viewpoint in the \nc setting  produces a simple set of criteria for a $2$-form, which when satisfied, gives a coherent system for constructing Hodge maps, metrics,  codifferentials, and Dirac operators. Moreover,  it also produces \nc generalisations of  classical K\"ahler phenomena that have not before appeared in the literature: Lefschetz decomposition, the Lefschetz identities,  Hodge decomposition, and the K\"ahler identities.

This is the first of a series of papers. In subsequent works we will enlarge the family of examples \cite{MM, SM}, investigate the analytic properties of some of the associated Dirac operators \cite{SM},  
and  investigate how the classical rules of Schubert calculus behave under $q$-deformation.  



The paper is organised as follows:  In Section 2 some well-known material is introduced about  Hopf algebras, quantum homogeneous spaces, and in particular quantum projective space. In Section 3 the theory  of covariant differential calculi  is recalled, as is the more recent notion of a complex structure.  The construction of the Heckenberger--Kolb calculus for $\cpn$ is also recalled, and some basic results about it presented.

In Section 4  symplectic  and  Hermitian structures are introduced. A \nc generalisation of Lefschetz decomposition is proved in {\bf Proposition \ref{LDecomp}} allowing for the definition of a Hodge map in {\bf Definition \ref{HDefn}}. A $q$-deformation of the fundamental form of the classical  Fubini--Study metric for complex projective space is then constructed, and a generalised version for the irreducible quantum flag manifolds is conjectured.
 
In Section 5 the construction of positive definite metrics from Hermitian forms is presented. Adjointability of $G$-comodule maps \wrt such metrics is then established in {\bf Corollary \ref{EqAdj}}, and presentations  of the codifferential and dual Lefschetz operators in terms of the Hodge map given. Finally, in  {\bf Proposition \ref{LIDS}} a deformed version of the Lefschetz identities is proved.

In Section 6, Hodge decomposition \wrt the holomorphic and anti-holomorphic derivatives is established, and shown to imply an isomorphism between cohomology classes and harmonic forms just as in the classical case. A \nc generalisation of Serre duality is also established.

In Section 7 the definition of a \nc K\"ahler structure is given and some of the basic results of classical geometry generalised, most notably the K\"ahler identities in {\bf Theorem \ref{KIDSII}}. Equality of the three Laplacians up to scalar multiple follows in {\bf Corollary \ref{LaplacianEq}}, implying in turn that  Dolbeault cohomology refines de Rham cohomology. Finally, a \nc generalisation of the hard Lefschetz theorem and the \linebreak $\del \adel$-lemma is given. The Hermitian form $\cpn$ is then observed to be K\"ahler, implying that the \hk calculus has cohomology groups of at least classical dimension.  We finish with some spectral calculations and a conjecture about constructing spectral triples from K\"ahler structures for the irreducible quantum flag manifolds.

Throughout the paper we endeavour to  present the derivation of all results as explicitly as possible, so as to make the paper accessible to a \ncg audience not necessarily familiar with  classical complex geometry.



\subsubsection*{Acknowledgements:} I would like to thank Karen Strung, Shahn Majid, Petr Somberg, Tomasz Brzezi\'nski, Edwin Beggs,  Vladimir Sou\v{c}ek,  Matthias Fischmann, and Adam-Christiaan van Roosmalen, for many useful discussions.

\section{Preliminaries on Quantum Homogeneous Spaces}

In this section we introduce some well-known material about cosemisimple Hopf algebras, quantum homogeneous spaces, and Takeuchi's categorical equivalence. The motivating example, quantum projective space, is also introduced.

 \subsection{Compact Quantum Group Algebras}

Let $G$ be a Hopf  \alg with comultiplication $\DEL$, counit $\e$, antipode $S$, unit $1$, and multiplication $m$. 
Throughout, we use Sweedler notation, and denote $g^+ := g - \e(g)1$, for $g \in G$, and $V^+ = V \cap \ker(\e)$, for $V$   a subspace of $G$. 

For any left $G$-comodule $(V,\DEL_L)$, its space of {\em matrix elements} is the co\alg
\bas
\C(V) : = \spn_\bC\{(\id \oby f)\DEL_L(v) \,|\, f \in \text{Lin}_{\bC}(V,\bC), v \in V\} \sseq G.
\eas
A comodule is irreducible \iff its co\alg of matrix elements is irreducible, and, for $W$ another left $G$-comodule,  $\C(V) = \C(W)$ \iff $V$ is equivalent to $W$. Moreover, $\C(V)$ decomposes as a left $G$-comodule into $\dim_{\bC}(V)$ copies of $V$. 

The notion of cosemisimplicity for a Hopf algebra will be essential in  the paper and all Hopf \algs  will be assumed to have the property. We present  three equivalent formulations of the definition (a proof of their equivalence can be found in \cite[\textsection 11.2.1]{KSLeabh}).

\begin{defn}
A Hopf algebra $G$ is called  {\em cosemisimple} if it satisfies the following three equivalent conditions:
\begin{enumerate}
\item It holds that $G \simeq \bigoplus_{V\in \wh{G}} \C(V)$, where summation is over all equivalence classes of left $G$-comodules.
\item Every comodule of $G$ is a direct sum of (necessarily finite) irreducible comodules.
\item  There exists a unique linear map $\haar:G \to \bC$, which we call the {\em Haar functional}, such that $\haar(1) = 1$, and 
\bas
(\id \oby \haar) \DEL(g) = \haar(g)1, & & (\haar \oby \id)\DEL(g) = \haar(g)1.
\eas
\end{enumerate}
\end{defn}

While the assumption of cosemisimplicity is enough for most of our requirements, we will need something stronger when discussing positive definiteness  in \textsection 5. 

\begin{defn}
A {\em compact quantum group algebra} is a cosemisimple Hopf $*$-algebra such that  $\haar(aa^*) > 0$, for all $a\neq 0$.
\end{defn}

The condition of  $G$ being a compact quantum group \alg is equivalent to it being the dense Hopf  \alg of a compact quantum group \cite{WOR2}, Woronowicz's celebrated structure  in the $C^*$-algebraic approach to quantum groups \cite{Wor}.

For any compact quantum group algebra, an inner product is given by the map
\bal \label{hinnerproduct}
G \oby G \to \bC, & & g \oby f \mto \haar(fg^*).
\eal
Moreover, \wrt this inner product, the decomposition  $G \simeq \bigoplus_{V \in \wh{G}} \C(V)$ is orthogonal.

\subsection{Quantum Homogeneous Spaces}

For a right $G$-comodule $V$ with coaction $\DEL_R$, we say that an element $v \in V$ is {\em coinvariant} if $\DEL_R(v) = v \oby 1$. We denote the subspace of all coinvariant elements by $V^G$, and call it the {\em coinvariant subspace} of the coaction. We also use the analogous conventions for left comodules. 
\begin{defn}
For $H$ a Hopf \algn, a {\em homogeneous} right $H$-coaction on $G$ is a coaction of the form $(\id \oby \pi)  \DEL$, where $\pi: G \to H$ is a surjective Hopf \alg map. A  {\em quantum homogeneous space}  $M:=G^H$ is the coinvariant subspace of such a coaction. 
\end{defn}

In this paper we will {\em always} use the symbols $G,H,\pi$ and $M$ in this sense. As is easily seen, 
 $M$ is a sub\alg  of $G$. Moreover, if $G$ and $H$ are Hopf $*$-\algn s, and $\pi$ is a Hopf $*$-\alg map, then $M$ is a $*$-sub\alg of $G$.

Our assumption of cosemisimplicity for Hopf algebras implies that $G$ is faithfully flat over $M$  \cite[Theorem 5.1.6]{PSHJS}. Recall that  $G$ is said to be {\em faithfully flat} as a right module over $M$ if  the functor  $G \oby_M -:\mm \to $ $\!_G \! {\mathrm{\, Mod}}$, from the category of left $M$-modules to the category of complex vector spaces, maps a sequence to an exact sequence if and only if the original sequence is exact. This is necessary in particular for the categorical equivalence of the next section to hold. 



\subsection{Takeuchi's Categorical Equivalence}

We now define the abelian categories $\lgmmm$ and $\lmhm$. The objects in $\lgmmm$ are  $M$-bimodules $\E$ (with left and right actions denoted by juxtaposition) endowed with a left $G$-coaction $\DEL_L$ \st  $\E M^+ \sseq M^+ \E$, and
\bal \label{gmmmdefn}
\DEL_L(m e m') =  m\1 e\m1 m'\1 \oby m\2 e\0 m'\2, & & \text{for all ~} m,m' \in M, e \in \E.
\eal
The morphisms in $\lgmmm$ are those  $M$-bimodule homomorphisms which are also homomorphisms of left $G$-comodules. The objects in $\lmhm$ are left $H$-comodules $V$  endowed with the trivial right $M$-action $(v,m) \mto \e(m)v$.  
The morphisms in $\lmhm$ are the left $H$-comodule maps. (Note that $\lmhm$ is equivalent under the obvious forgetful functor  to $^H\!\mathrm{Mod}$, the category of  left $H$-comodules.)   


If $\E \in \lgmmm$, then $\E/(M^+\E)$ becomes an object in $\lmhm$ with the left $H$-coaction 
\bal \label{comodstruc0}
\DEL_L[e] = \pi(e\m1) \oby [e\0], & & e \in \E,
\eal
where $[e]$ denotes the coset of $e$ in  $\E/(M^+\E)$. We define a functor $$\Phi:\lgmmm \to \lmhm$$  as follows:  
$\Phi(\E) :=  \E/(M^+\E)$, and if $g : \E \to \F$ is a morphism in $\lgmmm$, then $\Phi(g):\Phi(\E) \to \Phi(\F)$ is the map to which $g$ descends on $\Phi(\E)$.

If $V \in \lmhm$, then the {\em cotensor product} of $G$ and $V$, defined by
\bas
G \coby V := \ker(\DEL_R \oby \id - \id \oby \DEL_L: G\oby V \to G \oby H \oby V),
\eas  
becomes an object in $\lgmmm$ by defining an $M$-bimodule structure
\begin{align} \label{rightmaction}
m \Big(\sum_i g^i \oby v^i\Big)  = \sum_i m g^i \oby v^i, & & \Big(\sum_i g^i \oby v^i\Big) m = \sum_i g^i m \oby v^i,
\end{align}
and a left $G$-coaction 
\begin{align*}
\DEL_L\Big(\sum_i g^i \oby v^i\Big) = \sum_i g^i\1 \oby g^i\2 \oby v^i.
\end{align*}
We define a functor $\Psi:\lmhm \to \lgmmm$ as follows:
$
\Psi(V) := G \coby V,
$
and if $\g$ is a morphism in $\lmhm$, then $\Psi(\g) := \id \oby \g$. 

\begin{thm}\cite[Theorem 1]{Tak} An equivalence of the categories 
 $\lgmmm$ and  $\lmhm$, which we call {\em Takeuchi's equivalence},  is given by the functors $\Phi$ and $\Psi$ and the natural transformations
\begin{align}
\counit:\Phi \circ  \Psi(V) \to V, & & \Big[\sum_i g^i \oby v^i\Big] \mto \sum_{i} \e(g^i)v^i \label{counit},\\
\unit: \E \to \Psi \circ \Phi(\E), & & e \mto e\m1 \oby [e\0]. \label{unit}
\end{align}
\end{thm}

\begin{cor}
Takeuchi's equivalence restricts to an equivalence of categories between $\sgmm$  and  $\smh$, where $\sgmm$ is the full subcategory of $\lgmmm$ consisting of  finitely generated left $M$-modules, and $\smh$ is the full subcategory of $\lmhm$ consisting of  finite-dimensional comodules.

\end{cor}
\demo
We begin by recalling the well-known   \cite[\textsection 1]{Tak} isomorphism
\bas
G\oby_M \E \to G \oby \Phi(\E), & & g \oby_M e \mto ge\m1 \oby [e\0].
\eas 
This implies that, for any $V \in \smh$, we have that $G \oby_M \Psi(V)$ is finitely generated as a left $G$-module. 
Now cosemisimplicity of $H$ implies that there exists a projection  $\Phi(G) \to \Phi(M)$, and so, we have a projection $\r:G \to M$. The image of $G\oby_M \E$ under $m  (\r \oby \id)$ is isomorphic to $\E$ which we now see to be finitely generated. 
The proof of the converse is elementary.
\qed

We define the {\em dimension} of an object $\E \in \sgmm$ to be the vector space dimension of $\Phi(\E)$. Note that by cosemisimplicity of $G$, the abelian category $\smh$ is  semisimple, and so, $\sgmm$ is semisimple.

For $\E,\F$ two objects in $\sgmm$, we denote by $\E \oby_M \F$  the usual bimodule tensor product endowed with the standard left $G$-comodule structure.
It is easily checked that $\E \oby_M \F$ is again an object in $\sgmm$, and so, the tensor product $\oby_M$ gives the category a monoidal structure. With respect to the obvious monoidal structure on $\smh$, Takeuchi's equivalence is given the structure of a monoidal equivalence (see \cite[\textsection 4]{MMF2} for details) by the morphisms
\bas
\mu_{\E, \F}: \Phi(\E) \oby \Phi(\F) \to \Phi(\E \oby_M \F), & & [e] \oby [f] \mto [e \oby_M f], &  \text{ ~~~~~~~~ for } \E,\F \in \sgmm.
\eas
In what follows, this monoidal equivalence will be tacitly assumed.

Finally, we note that, for any $\E \in \sgmm$, the following decomposition exists:
\bas
\E \simeq G \coby \Phi(\E) \simeq \Big(\bigoplus_{V\in \wh{G}} \C(V))\Big) \coby \Phi(\E) = \bigoplus_{V\in \wh{G}} \big(\C(V) \coby  \Phi(\E)\big) =: \bigoplus_{V\in \wh{G}} \E_V.
\eas 
We call this the {\em Peter--Weyl decomposition} of $\E$.


\subsection{Quantum Projective Space}

We  recall the definition of the well-known quantum coordinate algebras $\bC_q[U_n]$ and $\bC_q[SU_n]$,   as well as the definition of quantum projective space, the motivating example considered throughout the paper. We finish with a discussion of weight decomposition for $\bC_q[U_n]$-comodules, an important tool in what follows.

\subsubsection{The Quantum Groups $\bC_q[U_n]$ and $\bC_q[SU_n]$}

We begin by fixing notation and recalling the various definitions and constructions needed to introduce the quantum unitary group and the quantum special unitary group. (Where proofs or basic details are omitted we refer the reader to \cite[\textsection 9.2]{KSLeabh}.)

For $q \in \bR_{>0}$, let $\bC_q[GL_n]$ be the quotient of the free \nc \alg  \linebreak ${\bC \big{<} u^i_j, \dt_n \inv \,|\, i,j = 1, \ldots , n \big{>}}$ by the ideal generated by the elements
\begin{align*}
u^i_ku^j_k  - qu^j_ku^i_k, &  & u^k_iu^k_j - qu^k_ju^k_i,                    & &   \; 1 \leq i<j \leq n, 1\leq k \leq n; \\
  u^i_lu^j_k - u^j_ku^i_l, &  & u^i_ku^j_l - u^j_lu^i_k - (q-q^{-1})u^i_lu^j_k, & &   \;  1 \leq i<j \leq n,\; 1 \leq k < l \leq n;\\
  \dt_n \dt_n \inv - 1,   & & \dt_n \inv \dt_n - 1,& &    
\end{align*}
where  $\dt_n$,  the {\em quantum determinant},  is the element
\[
\dt_{n} := \sum\nolimits_{\pi \in S_n}(-q)^{\ell(\pi)}u^1_{\pi(1)}u^2_{\pi(2)} \cdots u^n_{\pi(n)},
\]
with summation taken over all permutations $\pi$ of the set $\{1, \ldots, n\}$, and $\ell(\pi)$ is the number of inversions in $\pi$. As is well-known \cite[\textsection 9.2.2]{KSLeabh}, $\dt_n$ is a central element of the algebra.

A bi\alg structure on $\bC_q[GL_n]$ with coproduct $\DEL$, and counit $\e$, is uniquely determined by $\DEL(u^i_j) :=  \sum_{k=1}^n u^i_k \oby u^k_j$; $\DEL(\dt_n\inv) = \dt_n\inv \oby \dt_n\inv$; and $\e(u^i_j) := \d_{ij}$; $\e(\dt_n\inv) = 1$. The element $\det_n$ is grouplike \wrt $\DEL$ \cite[\textsection 9.2.2]{KSLeabh}. Moreover, we can endow  $\bC_q[GL_n]$ with a Hopf \alg structure by defining
\begin{align*}
S(\dt_{n} \inv) := \dt_{n}, ~~~~ S(u^i_j) := (-q)^{i-j}\sum\nolimits_{\pi \in S_{n-1}}(-q)^{\ell(\pi)}u^{k_1}_{\pi(l_1)}u^{k_2}_{\pi(l_2)} \cdots u^{k_{n-1}}_{\pi(l_{n-1})}\dt_n\inv,
\end{align*}
where $\{k_1, \ldots ,k_{n-1}\} := \{1, \ldots, n\}\bs \{j\}$, and $\{l_1, \ldots ,l_{n-1}\} := \{1, \ldots, n\}\bs \{i\}$ as ordered sets.
A Hopf $*$-\alg structure is determined by ${(\dt_{n}^{-1})^* = \dt_{n}}$, and $(u^i_j)^* =  S(u^j_i)$. We denote the Hopf $*$-\alg by $\bC_q[U_n]$, and call it the {\em quantum unitary group of order $n$}. We denote the Hopf  {$*$-\algn} $\bC_q[U_n]/\la \dt_{n} - 1 \ra$ by $\bC_q[SU_n]$, and call it the {\em quantum special unitary group of order $n$}.

\subsubsection{Quantum Projective Space}

Following the description introduced in \cite[\textsection 3]{Mey}, we present quantum  $n$-projective space as the  sub\alg of coinvariant elements of a $\bC_q[U_{n}]$-coaction on $\bC_q[SU_{n+1}]$.  (This sub\alg is a $q$-deformation of the coordinate \alg of the complex manifold $SU_{n+1}/ U_{n}$. Recall that  $\bC P^{n}$ is isomorphic to $SU_{n+1}/ U_{n}$.)

\begin{defn}
Let $\a_n:\bC_q[SU_{n+1}] \to \bC_q[U_{n}]$ be the surjective Hopf $\ast$-\alg map  defined by setting $\a_n(u^1_1) = \dt_{n} \inv,$ \, $\a_n(u^1_i)=\a_n(u^i_1)=0$, for $i = 2, \cdots, n+1$, and  $\a_n(u^i_j) = u^{i-1}_{j-1},$ for  $i,j=2, \ldots, n+1$. {\em Quantum projective $n$-space} $\cpn$ is defined to be the quantum homogeneous space of the corresponding homogeneous coaction \linebreak $(\id \otimes \a_n) \circ \Delta$. 
\end{defn}

As is well known  \cite[\textsection 11.6]{KSLeabh},  $\cpn$ is generated as a $\bC$-\alg by the set 
\bas
\{z_{ab} := u^a_1S(u^1_b)\,|\, a,b = 1,\ldots , n\}.
\eas 

\subsubsection{Weight Vectors for Objects in $\smh$}

Let $\bC[\bT^n]$ be the commutative polynomial \alg generated by  $t_k, t\inv_k$, for $k=1, \ldots, n$, satisfying the obvious relation $t_kt\inv_k  = 1$. We can give  $\bC[\bT^n]$ the structure of a Hopf  \alg by defining a coproduct, counit and antipode according to  $\DEL(t_k) := t_k \oby t_k$, $\e(t_k) := 1$, and $S(t_k) := t_k\inv$. Moreover, $\bC[\bT^n]$ has a Hopf $*$-\alg structure defined by $t^*_k := t_k\inv$. (Note that $\bC[\bT^1] \simeq \bC[U_1]$.)


A basis of $\bC[\bT^n]$ 
is given by 
\bas
\{t^\l := t^{l_1}_1 \cdots t^{l_n}_n \,|\, \l = (l_1, \ldots, l_n) \in \bZ^n\}.  
\eas
Since each basis element is grouplike, a $\bC[\bT^n]$-comodule structure is equivalent to a $\bZ^{n}$-grading. 
We  call the homogeneous elements of such a grading  {\em weight vectors}, and we call their degree their {\em weight}. 

We are interested in  $\bC[{\mathbb T}^n]$ because of the existence of the following map: Let  \linebreak $\t:\bC_q[U_{n}] \to \bC[\bT^{n}]$  be the surjective Hopf $\ast$-\alg map 
defined by 
\bas
\t(\dt_{n}\inv) := t_\bullet\inv, & & \t(u^i_j) := \d_{ij} \, t_i, & & \text{ for } i,j = 1, \dots, n,
\eas
where $t_{\bullet} :=  t_1 \cdots t_{n}$. For any left $\bC_q[U_{n}]$-comodule $V$, a left $\bC[\bT^{n}]$-comodule structure on $V$ is defined by $\DEL_{L,\t}:=$
$(\t \oby \id)  \DEL_L$.

\begin{lem}\label{weightsaddtens}
For any two objects $\D, \F \in \gmmm$, and  $d \in \D, f \in \F$ weight vectors of weight $w$ and $v$ respectively, then $d \oby_M f \in \D \oby_M \F$ is a weight vector of weight $w + v$.
\end{lem}
\demo
This follows directly from
\bas
~~~~~~~~~~ \DEL_L[e \oby_M f] = \t(d\m1f\m1) \oby [d\0 \oby_M f\0] = t^{w+v} \oby [d\0 \oby_M f\0].  ~~~~~~~~~~~~~ \text{\qed}
\eas

\section{Preliminaries and Basic Constructions on Differential Calculi and Complex Structures}

In this section we recall some well-known definitions from the theory of differential calculi, including material on $*$-calculi, orientability, and integrals. Some more recent material on complex structures is also considered. Finally, a concise presentation of the Heckenberger--Kolb calculus for $\cpn$ is given, and some elementary results on weight space decomposition proved.

\subsection{Complexes and Double Complexes}

For $(S,+)$ a commutative semigroup, an  {\em {$S$-graded \algn}} is an \alg of the form $\A = \bigoplus_{s \in S}\A^s$, where each $\A^s$ is a linear
subspace of $\A$, and $\A^s\A^t \sseq \A^{s+t}$, for all $s,t
\in S$. If $a \in \A^s$, then we say that $\a$ is a {\em homogeneous element of degree $s$}.  A {\em homogeneous mapping of degree $t$} on $A$ is a linear mapping $L:\A \to \A$ \st if $\a \in
\A^s$, then $L(\a) \in \A^{s+t}$. We say that a subspace $\B$ of $\A$ is {\em homogeneous} if it admits a decomposition  $\B = \oplus_{s \in S} \B^s$, with $\B^s \sseq \A^s$, for all $s \in S$.

A pair $(\A,\exd)$ is called a {\em complex} if $\A$ is an $\bN_0$-graded \algn, and $\exd$ is a homogeneous mapping  of degree $1$ \st $\exd^2 = 0$. A triple $(\A,\del,\ol{\del})$ is called a {\em double complex} if $\A$ is an $\bN^2_0$-graded \algn, $\del$ is homogeneous mapping of degree $(1,0)$, $\ol{\del}$ is homogeneous mapping of degree $(0,1)$, and 
\bas
\del^2 = \ol{\del}^2 = 0, & & \del  \ol{\del} =  - \ol{\del}  \del.
\eas 
Note we can  associate to any double complex $(\A,\del,\ol{\del})$ three different complexes  
\bas
(\A,\exd:=\del + \ol{\del}), & & (\A,\del), & & \text{ and  ~~~~~~ } & (\A,\ol{\del}),
\eas
where the $\bN_0$-grading on $\A$ is given by $\A^k := \bigoplus_{a+b=k} \A^{(a,b)}$. 

For any complex $(\A,\exd)$, we call an element  {\em $\exd$-closed} if it is contained in $\ker(\exd)$, and {\em $\exd$-exact} if it is contained in im$(\exd)$. Moreover, the {\em $\exd$-cohomology group} of order $k$ is the space
\bas
H^k_\exd := \ker(\exd:\A^k \to \A^{k+1})/\text{im}(\exd:\A^{k-1} \to \A^k).
\eas
For a double complex $(\A,\del, \ol{\del})$ we define $\del$-closed, $\adel$-closed, $\del$-exact, and $\adel$-exact forms analogously. The {\em$\del$-cohomology group} $H_\del^k$, and the {\em $\adel$-cohomology group} $H_{\adel}^k$,  are the cohomology groups of the complexes $(\A,\del)$ and $(\A,\ol{\del})$.  Finally, we note that we have the decompositions
\bas
& H_\del^k = \bigoplus_{a+b=k} H_{\del}^{(a,b)},& &  \text{ and } & &H_{\adel}^k = \bigoplus_{a+b=k} H_{\adel}^{(a,b)},&
\eas
where $H^{(a,b)}_\del$ and  $H^{(a,b)}_{\adel}$ are the $a\th$, and $b\th$,  cohomology groups of the complexes $(\A^{(\bullet,b)},\del)$ and $(\A^{(a,\bullet)},\adel)$ respectively, where the gradings are the obvious ones.

\subsubsection{Differential $*$-Calculi}

A complex $(\A,\exd)$ is called a {\em differential graded algebra} if $\exd$ is a {\em graded derivation}, which is to say, if it satisfies the {\em graded Leibniz rule}
\bas
\exd (\a\b) = \exd (\a)\b+(-1)^k \a \exd \b,       &  & \textrm{  for all $\a \in \A^k$, $\b \in \A$}.
\eas
The operator $\exd$ is called the {\em differential} of the differential graded \algn. 
\begin{defn}
A {\em differential calculus} over an \alg $A$ is a differential graded
\alg  $(\Om^\bullet,\exd)$ \st $\Om^0=A$, and
\begin{align} \label{pofu}
\Om^k = \spn_{\bC}\{a_0\exd a_1\wed  \cdots \wed \exd a_k \,|\, a_0, \ldots, a_k \in A\}.
\end{align}
\end{defn}
We use $\wedge$ to denote the multiplication between elements of a differential calculus when both are of order greater than  $0$. We call an element of a differential calculus a {\em form}. 
A {\em differential map} between two differential calculi $(\Om^\bullet,\exd_{\Om})$ and $(\G^\bullet,\exd_{\G})$, defined over the same \alg $A$, is a bimodule map $\f: \Om^\bullet  \to \Gamma^\bullet$ \st $\f \circ \exd_{\Om} = \exd_{\Gamma}$. 

\bigskip

We call a differential calculus $(\Om^\bullet, \exd)$ over a $*$-\alg $A$ a {\em  $*$-differential calculus} if the involution of $A$ extends to an involutive conjugate-linear map on $\Om^\bullet$, for which $(\exd \w)^* = \exd \w^*$, for all $\w \in \Om$, and
\[
(\w \wed \nu)^*=(-1)^{kl}\nu^* \wed \w^*, \qquad \text{ for all }\w \in
\Om^k,~ \nu \in \Om^l.
\]
We say that a form $\w \in \Om^\bullet$ is {\em real}  if $\w^* = \w$.

\bigskip

A differential calculus $\Omega^\bullet$ over a \qhs $M$ is said to be {\em covariant}  if $\DEL_L: M \to G \oby M$ extends to a necessarily unique \alg map  $\DEL_L: \Om^\bullet \to G \oby \Om^\bullet$ \st 
\bas
\DEL_L(m\exd n) = \DEL_L(m) \big((\id \oby \exd) \circ \DEL_L(n)\big) = m\1 n\1 \oby m\2 \exd n\2,  &  & m,n \in M.
\eas 
In this paper, all covariant calculi will be assumed to be finite-dimensional and to satisfy $\Om^\bullet M^+ \sseq M^+\Om^\bullet$, giving  $\Om^\bullet$ the structure of an object in $\sgmm$. This implies that a multiplication is defined on $\Phi(\Om^\bullet)$  by $[\w] \wed [\nu] := [\w \wed \nu]$. It follows from (\ref{pofu}) that every element of $\Phi(\Om^k)$ is a sum of elements of the form $[\w_1] \wed \cdots \wed [\w_k]$, for $w_i \in \Om^1$. When working with covariant calculi we usually use the convenient notation $V^\bullet := \Phi(\Om^\bullet)$.

\subsection{Orientability and Closed Integrals}

We say that a differential calculus has {\em total dimension} $n$ if $\Om^k = 0$, for all $k>n$, and $\Om^n \neq 0$. If in addition there exists an $A$-$A$-bimodule isomorphism $\vol: \Om^n \simeq A$,  then we say that $\Om^{\bullet}$ is {\em orientable}. We call a choice of such an isomorphism an {\em orientation}.  If $\Om^\bullet$ is a covariant calculus over a quantum homogeneous space $M$ and $\vol$ is a morphism in $\sgmm$, then we say that $\Om^\bullet$ is  {\em covariantly orientable}. Note all covariant orientations  are equivalent up to  scalar multiple. If  $\Om^\bullet$ is a $*$-calculus over a $*$-\algn, then a \mbox{{\em $*$-orientation}} is an orientation which is also a $*$-map. A {\em $*$-orientable calculus} is one which  admits a $*$-orientation.

When the calculus is defined over a \qhsn, we define the {\em integral}, \wrt $\vol$, to be the map which is zero on all $\Om^k$, for $k<n$, and
\bas
\int: \Om^n \to \bC, & & \w \mto \haar(\vol(\w)),
\eas
where $\haar$ is the Haar functional. We say that the integral is {\em closed} if 
$
\int \exd \w = 0$,  \text{for all } $\w \in \Om^{n-1}.
$


\begin{lem}\label{RiemClosedInt}
For $\Om^\bullet$ a covariant orientable calculus over a quantum homogeneous space $M$, the integral is closed \iff  $\exd\big(\,^G\big(\Om^{n-1}\big)\big) = 0$.
\end{lem}
\demo
Cosemisimplicity of $G$ guarantees that we can make a choice of complement $K \in \, ^G$Mod to $^G\big(\Om^{n-1}\big)$ in $\Om^{n-1}$. Since the map $\int \circ \, \exd: \Om^{n-1} \to \bC$ is a left $G$-comodule map, its restriction to $K$ must be the zero map.  Hence if $\exd(^G\Om^{n-1}) = 0$, then  $\int \exd \w = 0$, for all $\w \in \Om^{n-1}$. 

Conversely, for any $\w \in \,^G\big(\Om^{n-1}\big)$, the fact that $\exd$ is a comodule map implies that $\exd \w = \l \vol\inv(1)$, for some $\l \in \bC$. Moreover, $\int \exd \w = \int \l \vol\inv(1) = \l \haar(1) = \l$. Hence, if $\l \neq 0$, which is to say if $\exd \w \neq 0$, then the integral is not closed.
\qed

\begin{cor} \label{closurev}
The integral is closed if the decomposition of $V^{2n-1}$ into irreducible comodules does not contain the trivial comodule.
\end{cor}
\demo
If there exists a left $G$-coinvariant form $\w \in \Om^{n-1}$, then $\Phi(M\w)$ is a trivial subcomodule of $V^{n-1}$. Hence, if no such subcomodule exists, there can be no coinvariant forms, and $d\big(\,^G\Om^{(n-1)}\big) =0$ is satisfied vacuously.
\qed

\subsection{Complex Structures}

In this subsection we recall the basic definitions and results of complex structures. For a more detailed introduction see \cite{MMF2}.

\begin{defn}\label{defnnccs}
An {\em almost complex structure} for a  differential $*$-calculus  $\Om^{\bullet}$, over a \mbox{$*$-\alg $A$}, is an $\bN^2_0$-\alg grading $\bigoplus_{(a,b)\in \bN^2_0} \Om^{(a,b)}$ for $\Om^{\bullet}$ such that 
\begin{enumerate}
\item \label{compt-grading}  $\Om^k = \bigoplus_{a+b = k} \Om^{(a,b)}$, ~~~~ for all $k \in \bN_0$,
\item  \label{star-cond} $(\Om^{(a,b)})^* = \Om^{(b,a)}$, ~~~~~~~~ for all $(a,b) \in \bN^2_0$.
\end{enumerate}
\end{defn}
We call an element of $\Om^{(a,b)}$ an {\em  $(a,b)$-form}. We say that an almost complex structure is {\em factorisable} if we have bimodule isomorphisms
\bal  \label{wedge-cond} 
\wed:\Om^{(a,0)} \oby_A \Om^{(0,b)} \simeq \Om^{(a,b)}  & & \text{ and } & & \wed: \Om^{(0,b)} \oby_A \Om^{(a,0)} \simeq \Om^{(a,b)}.
\eal If the \alg is a \qhs and $\Om^\bullet$ is a covariant calculus,  then we say that the almost complex structure is {\em covariant} if $\Om^{(a,b)}$ is a sub-object of $\Om^\bullet$ in $\sgmm$, for all  $(a,b) \in \bN^2_0$. We say that an almost complex structure is of {\em diamond type}  if whenever  $a > n$, or $b > n$, then necessarily  $\Om^{(a,b)} = 0$. Note that any almost complex structure on the de Rham complex of a manifold is necessarily of diamond type. Moreover, any calculus which is not of diamond type can clearly be quotiented to a calculus of diamond type.

Let $\del$ and $\ol{\del}$ be the unique homogeneous operators  of order $(1,0)$, and $(0,1)$ respectively, defined by 
\bas
\del|_{\Om^{(a,b)}} : = \proj_{\Om^{(a+1,b)}} \circ \exd, & & \ol{\del}|_{\Om^{(a,b)}} : = \proj_{\Om^{(a,b+1)}} \circ \exd,
\eas
where $ \proj_{\Om^{(a+1,b)}}$, and $ \proj_{\Om^{(a,b+1)}}$, are the projections from $\Om^{a+b+1}$ onto $\Om^{(a+1,b)}$, and $\Om^{(a,b+1)}$, respectively.  Assuming that the calculus is of total dimension $2n$, and that the almost complex structure is of diamond type, then $\exd$ restricts to $\del$ on $\Om^{(n-1,n)}$, and to $\adel$ on $\Om^{(n,n-1)}$. Hence closure of the integral is equivalent to 
\bal \label{complexcloseure}
\int \del \w = \int \adel \w' = 0, & & \text{ for all }  \w \in \Om^{(n-1,n)}, \, \w' \in \Om^{(n,n-1)}. 
\eal

As observed in \cite[\textsection 3.1]{BS} the proof of the following lemma carries over directly from the classical setting \cite[\textsection 2.6]{HUY}. Since the formulation of the definition of an almost complex structure used here differs from that used in \cite[\textsection 3.1]{BS} (see Remark \ref{rembs} below) we include a proof.

\begin{lem}\label{intlem} \cite[\textsection 3.1]{BS}
If $\bigoplus_{(a,b)\in \bN^2_0}\Om^{(a,b)}$ is an almost complex structure for a differential $*$-calculus $\Om^{\bullet}$ over an \alg $A$, then the following two conditions are equivalent:
\begin{enumerate} 
\item $\exd = \del + \ol{\del}$,
\item the triple $\big(\bigoplus_{(a,b)\in \bN^2}\Om^{(a,b)}, \del,\ol{\del}\big)$ is a double complex.
\end{enumerate}
\end{lem}
\demo
Let us first show that {\em 1} implies {\em 2}.  For any $\w \in \Om^k$,
\begin{align*}
0 = \exd^2(\w) =  \del^2(\w) + (\ol{\del} \circ \del +
\del \circ  \ol{\del})(\w) + \ol{\del}^2(\w).
\end{align*}
Since each of the three summands on the right hand side lie in complementary
subspaces of $\Om^{k+2}(M)$,  each must be zero.

Let us now show that {\em 2} implies {\em 1}. Note first that, for $g \in A$, 
\bas
0 = \proj_{\Om^{(0,2)}}(\exd^2 g) = & \, \proj_{\Om^{(0,2)}}\big(\exd(\del g + \adel g)\big) \\
 = & \, \proj_{\Om^{(0,2)}}\big(\exd(\del g)\big) + \adel^2 g \\
= & \, \proj_{\Om^{(0,2)}}(\exd(\del g)). 
\eas
Thus, for any $f \in A$, the form $\exd(f \del g) = \exd f \wed \del g + f \exd(\del g)$  is contained in  $\Om^{(2,0)} \oplus \Om^{(1,1)}$, and so,
$
\exd(\Om\hol) \sseq \Om^{(2,0)} \oplus \Om^{(1,1)}.
$
An analogous argument, using instead the projection $\proj_{\Om^{(2,0)}}$, shows that $
\exd(\Om\hol) \sseq \Om^{(2,0)} \oplus \Om^{(1,1)}.
$
Since $\Om^{(a,b)}$ is spanned by products of $a$ elements of $\Om^{(1,0)}$, and $b$ elements of $\Om^{(0,1)}$, it  follows from the Leibniz rule that $\exd = \del + \adel$ as required. 
\qed

\begin{defn}
When the conditions in Lemma \ref{intlem} hold for an almost complex structure, then we say that it is {\em integrable}. 
\end{defn}

We usually call an integrable almost complex structure a {\em complex structure}, and the double complex  $(\bigoplus_{(a,b)\in \bN^2}\Om^{(a,b)}, \del,\ol{\del})$  its {\em Dolbeault double complex}. An easy consequence of integrability is that 
\bal \label{stardel}
\del(\w^*) = (\adel \w)^*, & &  \adel(\w^*) = (\del \w)^*, & & \text{ for all  } \w \in \Om^\bullet.
\eal

\begin{rem}
The property of integrability for an almost complex structure has a number of other equivalent formulations in addition to the two presented above. For details see \cite[Lemma 2.13]{MMF2}. 
\end{rem}

\begin{rem}\label{rembs}
For a discussion of the equivalence of the definition of an almost complex structure used here with the definition of Beggs and Smith in \cite[Definition 2.6]{BS} see \cite[Remark 2.16]{MMF2}. 
\end{rem}

\subsection{The \hk Calculi for Quantum Projective Space}

In this subsection we give a brief presentation of the Heckenberger--Kolb calculus over $\cpn$. A more detailed presentation, in the notation of this paper, can be found in \cite{MMF2}. The calculi were originally introduced by Heckenberger and Kolb in $\cite{HK}$ for the more general class of examples given by the irreducible quantum flag manifolds, as discussed in \textsection 4.5. Their maximal prolongations and  complex structures were first presented in $\cite{HKdR}$.


\subsubsection{First-Order Calculi and Maximal Prolongations}

In this subsection we recall some details about first-order differential calculi necessary for our presentation of the \hk calculus below. A {\em first-order differential calculus} over $A$ is a pair $(\Om^1,\exd)$, where $\Omega^1$ is an $A$-$A$-bimodule and $\exd: A \to \Omega^1$ is a linear map for which the {\em Leibniz rule}, 
$
\exd(ab)=a(\exd b)+(\exd a)b, \text{ for } a,b,\in A,
$
holds and for which $\Om^1 = \spn_{\bC}\{a\exd b\, | \, a,b \in A\}$.  The notions of differential map,  and left-covariance when the calculus is defined over a \qhs $M$, have obvious first-order analogues, for details see \cite[\textsection 2.4]{MMF2}. The {\em direct sum} of two first-order differential calculi $(\Om^1,\exd_{\Om})$ and $(\G^1,\exd_{\G})$ is the first-order calculus $(\Om^1 \oplus \Gamma^1, \exd_\Om + \exd_\G)$.  Finally, we say that  a left-covariant first-order calculus over $M$ is {\em irreducible} if it does not possess any non-trivial quotients by a left-covariant $M$-bimodule. 

We say that a differential calculus $(\G^\bullet,\exd_\G)$ {\em extends} a first-order calculus $(\Om^1,\exd_{\Om})$ if there exists a bimodule isomorphism $\f:\Om^1 \to \G^1$ \st $\exd_\G = \f \circ \exd_{\Om}$. It can be shown  \cite[\textsection 2.5]{MMF2} that any first-order calculus admits an extension $\Om^\bullet$ which is maximal  in the sense that there exists a unique differential map from $\Om^\bullet$ onto any other extension of $\Om^1$. We call this extension the {\em maximal prolongation} of the first-order calculus.

\subsubsection{Defining the Calculus}

We present the calculus in two steps, beginning with Heckenberger and Kolb's classification of first-order calculi over $\cpn$, and then discussing the maximal prolongation of the direct sum of the two calculi identified. 

\begin{thm} \cite[\textsection 2]{HK} \label{HKClass}
There exist exactly two non-isomorphic irreducible left-covariant first-order differential calculi of finite dimension over $\cpn$. We call the direct sum of these two calculi the {\em Heckenberger--Kolb} calculus of $\cpn$.
\end{thm}

We denote these two calculi by $\Om\hol$ and $\Om\ahol$, and denote their direct sum by $\Om^1$. For a proof of the following lemma see \cite[Lemma 5.2]{MMF1}. 

\begin{lem}\label{basislem}
A basis of $V\hol:=\Phi(\Om\hol)$ and $V\ahol:=\Phi(\Om\ahol)$ is given respectively by 
\bas
\{ e^+_a := [\del z_{a+1,1}] \,  \,|\, a=1,\ldots, n \}, & &  \{e^-_a := q^{2(a+1)} [\adel z_{1,a+1}]\ \,|\, a = 1,\ldots, n\}.
\eas
Moreover, $[\del z_{ab}] = [\adel z_{ab}] = 0$ for $a,b \geq 2$ or $a = b = 1$.
\end{lem}

We call a subset $\{i_1, \ldots, i_k\} \sseq \{1,\ldots, n\}$ {\em ordered} if $a_1 < \cdots < a_k$. For any two ordered subsets $I, J \sseq  \{1,\ldots, n\}$, we denote 
\bas
e^+_I \wed e^-_J := e^+_{i_1} \wed \cdots \wed e^+_{i_k} \wed  e^-_{j_1} \wed \cdots \wed e^-_{j_{k'}}.
\eas

\begin{lem}\label{CalcBasis} For $\Om^\bullet$ the \mpr of $\Om^1$, the space  $V^{k}:= \Phi(\Om^k)$ has dimension  $\binom{2n}{k}$. A basis is given by 
\begin{align*}
\{ e^-_I \wed e^+_J \, | \, I,J \sseq \{1,\ldots, n\} \text{  ordered subsets  \st } |I|+|J| = k \}.
\end{align*}
\end{lem}

A full set of generating relations of the \alg $V^\bullet$ is given  in following lemma, for a proof see \cite[Proposition 5.8]{MMF2}.

\begin{prop}\label{prop:i2:HK}
The set of relations in $V^\bullet$ is generated by the elements
\bas
e^-_i \wed e^+_j + q  e^+_j \wed e^-_i, & & e^-_{i}\wed e^+_{i} + q^{2} e^+_{i}  \wed e^-_{i} +  (q^{2}-1) \sum_{a=i+1}^{n}  e^+_{a} \wed e^-_{a}, \\
e^-_i \wed e^-_h + q \inv e^-_h \wed e^-_i, & & e^+_i \wed e^+_h + q  e^+_h \wed e^+_i, ~~~~~ 
e^+_i \wed e^+_i, ~~~~ e^-_i \wed e^-_i, 
\eas
for  $h,i,j = 1, \ldots, n$, $i \neq j$, and $h< i$. 
\end{prop}


\subsubsection{Weight Space Decomposition of $V^1$}    

In this subsection we give an explicit description of the left module structure of $\Phi(\Om^1)$, as well as its weight space decomposition.  A proof of the $\bC_q[U_{n}]$-coaction formulae can be found in \cite[Lemma 6.11]{MMF2}. The weight decomposition is novel, and so, we include a proof.

\begin{lem}
The left $\bC_q[U_{n}]$-coactions  on  $V^{(1,0)}$ and $V^{(0,1)}$ are given by
\bas
\DEL_{L}(e^+_{i}) =  \sum_{k=1}^{n} u^{i}_{k}\dt_{n} \oby \,e^+_{k},& & \DEL_L(e^-_{i}) =  \sum_{k=1}^{n} S(u^{k}_{i})\dt_{n} \inv \oby \, e^-_{k}. 
\eas
\end{lem}

%

\begin{cor}\label{firstorder2}
The induced left $\bC_q[\bT^{n}]$-coactions  on  $V^{(1,0)}$ and $V^{(0,1)}$ are given by
\bas
\DEL_{L,\t}(e^+_{i}) =   t_i t_{\bullet} \oby e^+_{i},& & \DEL_{L,\t}(e^-_{i}) =    (t_i t_{\bullet}) \inv \oby e^-_{i}. 
\eas
\end{cor}
\demo
The first identity follows immediately from
\bas
\DEL_{L,\t}(e^+_{i}) =  & \sum_{k=1}^{n} \t(u^{i}_{k} \dt_{n}) \oby e^+_{k} =  \t(u^{i}_{i})t_{\bullet} \oby e^+_{i}  = t_i t_{\bullet} \oby e^+_i.
\eas
The second identity is established similarly.
\qed


\subsubsection{A Complex Structure}

Finally, we come to the definition of a complex structure for the calculus. Denote 
\begin{align*}
V^{(a,b)} := \spn_\bC\{ e^+_I \wed e^-_J \, | \, I,J \sseq \{1, \ldots, n\} \text{ ordered subsets with } |I|=a,\, |J| = b \}.
\end{align*}
The decomposition $V^k$ $\simeq \bigoplus_{(a+b = k)} V^{(a,b)}$, for all $k$, follows immediately. For a proof of the following proposition see \cite[\textsection 6, \textsection7]{MMF2}.

\begin{prop} 
For the Heckenberger--Kolb calculus over $\cpn$, there is a unique covariant complex structure  $\Om^{(\bullet,\bullet)}$ \st 
$
\Phi(\Om^{(a,b)}) = V^{(a,b)}.
$
Moreover, the complex structure is factorisable.
\end{prop}

\begin{lem}\label{multformweights} 
Every zero weight vector of $V^\bullet$ is contained in $\bigoplus_{k =1}^n V^{(k,k)}$.
\end{lem}
\demo
The lemma follows from Corollary \ref{firstorder2} and the multiplicativity of $\DEL_{L,\t}$. 
\qed

\begin{cor}
The integral associated to any covariant orientation of $\Om^\bullet$ is closed.
\end{cor}
\demo
The lemma tells us that $V^{(n-1,n)}$ and $V^{(n,n-1)}$ contain no zero weights. Hence, they contain no elements coinvariant \wrt $\DEL_L$. Closure of the integral now follows from Corollary \ref{closurev}.
\qed


\section{Hermitian Structures and Hodge Maps}

In this section we introduce  symplectic and Hermitian forms and use them to prove a \nc generalisation of the Lefschetz decomposition. Motivated by Weil's well-known classical  formula \cite[\textsection 1.2]{Weil} relating the Hodge map with Lefschetz decomposition, we  introduce a definition for a Hodge map associated to any Hermitian form. 

Throughout this section $\Om^\bullet$ denotes a differential calculus, over an \alg $A$, of total dimension $2n$.

\subsection{Almost Symplectic Forms}

As a first step towards the definition  of an Hermitian form, we present a direct \nc generalisation of the classical definition of an almost symplectic form. 

\begin{defn}
An {\em almost symplectic form} for $\Om^\bullet$ is a central real $2$-form $\s$ such that, \wrt the {\em Lefschetz operator}
\bas
L:\Om^\bullet \to \Om^\bullet,  & &   \w \mto \s \wed \w,
\eas
isomorphisms are given by
\bal \label{Liso}
L^{n-k}: \Om^{k} \to  \Om^{2n-k}, & & \text{ for all } 1 \leq k < n.
\eal
\end{defn} 

Note that since $\s$ is a central real form, $L$ is an $A$-$A$-bimodule $*$-homomorphism. Moreover, if $\s$ is an almost symplectic form for a covariant calculus over a quantum homogeneous space $M$, then $L$ is a morphism in $\sgmm$ \iff $\s$ is a left $G$-coinvariant form.

The existence of a symplectic form has important implications for the structure of a differential calculus. Crucial to understanding this structure is the notion of a primitive form, which directly generalises the classical definition of a primitive form \cite[\textsection 1.2, \textsection 3.1]{HUY}.

\begin{defn}\label{defn:sympform}
For $L$ the Lefschetz operator of any almost symplectic form,  the space of {\em primitive $k$-forms}  is
\bas
P^k : = \{\a \in \Om^{k} \,|\, L^{n-k+1}(\a) = 0\},  \text{ ~ if } k \leq n,& & \text{ and } & & P^k := 0, \text{ ~ if } k>n.
\eas
\end{defn}

One of the main reasons primitive forms are so important is the following decomposition result. It shows that  an almost symplectic form implies the existence  of a further refinement of the $\bN_0$-decomposition of a differential calculus.

\begin{prop}\label{LDecomp}
For $L$ the Lefschetz operator of any almost symplectic form, we have the  $A$-bimodule decomposition
\bas
\Om^k \simeq \bigoplus_{j \geq 0} L^j(P^{k-2j}),
\eas
which we call the {\em Lefschetz decomposition.}
\end{prop}
\demo
Let us assume that the decomposition holds for some $k \leq n-2$. Consider the composition  
\begin{displaymath}
    \xymatrix{
       \Om^{k}\ar@/^2.3pc/[rrrr]^{L^{n-k}} \ar[rr]_{L}
      &  &  \Om^{k+2} \ar[rr]_{L^{n-k-1}} &  & \Om^{2n-k}.
        }
\end{displaymath}
Since $L^{n-k}:\Om^{k} \to \Om^{2n-k}$ is an isomorphism of $A$-$A$-bimodules, we have the $A$-$A$-bimodule decomposition 
\bas
\Om^{k+2} \simeq & \, \ker\big(L^{n-k-1}|_{\Om^{k+2}}\big) \oplus  L(\Om^{k}) = \ker\big(L^{n-(k+2)+1}|_{\Om^{k+2}}\big) \oplus   L(\Om^{k})\\
= & \, P^{k+2} \oplus L(\Om^{k}) =  P^{k+2} \oplus \Big(\bigoplus_{j \geq 0} L^{j+1}(P^{k-2j})\Big) = P^{k+2} \oplus \Big(\bigoplus_{j \geq 1} L^{j}(P^{k+2-2j})\Big)\\
= & \bigoplus_{j \geq 0} L^j(P^{k + 2 - 2j}).
\eas
Since $\Om^0 = P^0$ and $\Om^1 = P^1$, it now follows from an inductive argument that the proposition holds for each space of forms of degree less than or equal to $n$.

Turning to forms of degree greater than $n$, we see that, for $k=0,\cdots, n$,
\bas
 \Om^{2n-k} \simeq & \, L^{n-k}(\Om^{k}) \simeq  L^{n-k}\Big(\bigoplus_{j \geq 0} L^j(P^{k-2j})\Big) = \bigoplus_{j \geq n-k} L^{j}(P^{2n-k-2j})  \\
= &  \bigoplus_{j \geq 0} L^{j}(P^{2n-k-2j}), 
\eas
where the last equality follows from the fact that, for $j=0,\dots, n-k-1$, either $2n-k-2j > n$ and $P^{2n-k-2j} = 0$ by definition, or $k+2 \leq 2n-k-2j \leq n$, and so, we have $L^j(P^{2n-k-2j}) = 0$. \qed

\subsection{Closed, Central,  and Symplectic Forms}

In general, it can prove tedious to verify that a given $2$-form is central. Assuming that the form is $\exd$-closed, however, makes the task much easier.

\begin{lem}\label{closcent}
A $\exd$-closed form  is central \iff it commutes with $0$-forms. 
\end{lem}
\demo
If $\s$ is a $\exd$-closed form which commutes with $0$-forms, then 
\bas
\s \wed \exd a = \exd(\s a) = \exd(a\s) = \exd a  \wed \s, & & \text{ for all~~}  a \in \Om^0.
\eas
Thus $\s$ commutes with all $1$-forms, and hence with all forms. The proof in the other direction is trivial.
\qed


\begin{lem}
If $\Om^\bullet$ is a covariant calculus over a quantum homogeneous space $M$,  then every left $G$-coinvariant form commutes with $0$-forms.
\end{lem}
\demo
With respect to the isomorphism $\unit:\Om^\bullet \simeq G \coby V^\bullet$, any left $G$-coinvariant  $\w$ satisfies  $\unit(\w) = 1 \oby [\w]$. That $m \in M$ commutes with $\w$ is obvious from (\ref{rightmaction}).
\qed

\begin{cor} \label{closedeqcentral}
Every left $G$-coinvariant $\exd$-closed form is central. 
\end{cor}

The following lemma gives us a sufficient  criterion  for  a coinvariant form to be $\exd$-closed. It should be noted, however, that the condition is not necessary.

\begin{lem} \label{v3closedform}
If $^H\!(V^{3})$ is trivial, then every left $G$-coinvariant $2$-form is $\exd$-closed.
\end{lem}
\demo
For a left $G$-coinvariant $2$-form $\w$, covariance of the calculus implies that $\DEL_L(\exd \w) = (\id \oby \exd) \DEL_L(\w) = 1 \oby \exd \w$.  Hence, if $\exd \w \neq 0$, the space  $^G(\Om^3)$ contains a non-trivial left $G$-coinvariant element. Since this cannot happen if $^H\!(V^{3})$ is trivial, we must conclude that $\exd \w = 0$. 
\qed

Motivated by this discussion of closed forms, we present the following \nc generalisation of the classical notion of a symplectic form \cite[\textsection 3.1]{HUY}.

\begin{defn}
A {\em symplectic form} is a $\exd$-closed almost symplectic form.
\end{defn}

As we will see in \textsection $7$, a K\"ahler form is a special type of symplectic form whose existence has many far-reaching consequences for the structure of a differential calculus.


\subsection{Hermitian Structures and  $h$-Hodge Maps}

We begin by introducing the notion of an Hermitian structure for a differential $*$-calculus, which is  essentially just a symplectic form interacting with a complex structure in a natural way.  In the commutative case each such form is the fundamental form of a uniquely identified Hermitian metric  \cite[\textsection 3.1]{HUY}.

\begin{defn} An {\em Hermitian structure} for a $*$-calculus $\Om^{\bullet}$  is a pair $(\Om^{(\bullet,\bullet)}, \s)$ where $\Om^{(\bullet,\bullet)}$ is a complex structure  and  $\s$ is an almost symplectic form, called the {\em Hermitian form}, \st $\s \in \Om^{(1,1)}$.
\end{defn}

When $\Om^\bullet$ is a covariant $*$-calculus over a \qhsn, $\Om^{(\bullet,\bullet)}$ is a covariant complex structure, and $\s$ is a left $G$-coinvariant form, then we say that $(\Om^{(\bullet,\bullet)}, \s)$  is a {\em covariant} Hermitian structure. We omit the proof of the following lemma which is clear.

 \begin{lem}\label{diamond}
The existence of an Hermitian structure for a complex structure implies that it is of diamond type.
\end{lem}

Taking our motivation from Weil's well-known classical formula \cite[\textsection 1.2]{Weil} presenting the Hodge map in terms of the Lefschetz decomposition, we use Lemma \ref{LDecomp} to introduce a \nc generalisation of the Hodge map. (Note that quantum integers are used instead of integers, as is discussed in the remark below.)

\begin{defn} \label{HDefn}
For $h \in \bR_{> 0}$, the {\em $h$-Hodge map} associated to an Hermitian structure is the morphism uniquely defined by
\bas
\ast_{h}(L^j(\w)) = (-1)^{\frac{k(k+1)}{2}}i^{a-b}\frac{[j]_h!}{[n-j-k]_h!}L^{n-j-k}(\w), & & \w \in P^{(a,b)} \sseq P^k,
\eas
where  $[m]_h := h^{m-1} + h^{m-3} + \cdots + h^{-m+1}$ denotes the quantum integer of $m$. We call $h$ the {\em Hodge parameter} of the Hodge map.
\end{defn}

As a first consequence of the definition, we establish direct generalisations of four of the basic properties of the classical Hodge map.

\begin{lem}\label{Hodgeprop} It holds that
\begin{enumerate}

\item $\ast_h^2(\w) = (-1)^k\id$, \text{ for all } $\w \in \Om^k$,

\item $\ast_h$ is an isomorphism,

\item $\ast_h(\Om^{(a,b)}) = \Om^{(n-b,n-a)}$,  
 \label{astpq}

\item $\ast_h$ is a $*$-map.

\end{enumerate}
\end{lem}
\demo
\begin{enumerate}

\item  By  Lefschetz decomposition it suffices to prove the result for a form of type $L^j(\a)$, for $\a \in P^{(a,b)} \sseq P^k$, $j \geq 0$. From the definition of $\ast_h$ we have that
\bas 
& \ast^2_h(L^j(\a))\\
 = & (-1)^{\frac{k(k+1)}{2}}i^{a-b}\frac{[j]_h!}{[n-j-k]_h!}\ast_h\big(L^{n-j-k}(\a)\big) \\
                                = & (-1)^{\frac{k(k+1)}{2}}i^{2(a-b)} \frac{[j]_h!}{[n-j-k]_h!}(-1)^{\frac{k(k+1)}{2}}\frac{[n-j-k]_h!}{[n-(n-j-k)-k]_h!}L^{n-(n-j-k)-k}(\a)\\
= & (-1)^k\frac{[j]_h!}{[j]_h!}L^j(\a) =  (-1)^k L^j(\a).
\eas

\item This follows immediately from 1.

\item This follows directly from the definition of $\ast_h$, the fact that $L$ is a degree $(1,1)$ map, and, as just established, the fact that $\ast_h$ is an isomorphism.

\item Again, it suffices to prove the result for a form of type $L^j(\a)$. Since $(P^{(a,b)})^* = P^{(b,a)}$, the definition of $\ast_h$ implies that
\bas
\ast_h  \big((L^j(\a))^*\big) & = \ast_h  (L^j(\a^*)) = (-1)^{\frac{k(k+1)}{2}}\frac{[j]_h!}{[n-j-k]_h!} i^{b-a} \big(L^{n-j-k}(\a^*)\big)\\
 & = \Big((-1)^{\frac{k(k+1)}{2}}\frac{[j]_h!}{[n-j-k]_h!} i^{a-b} L^{n-j-k}(\a)\Big)^*\\
& =  \big( \ast_h L^j(\a)\big)^*~~~~~~~~~~~~~~~~~~~~~~~~~~~~~~~~~~~~~~~~~~~~~~~~~~~~~~~~~~~~~~~~~~~ \text{\qed}
\eas

\end{enumerate}

\begin{cor}
With respect to a choice of Hermitian structure, a $*$-orientation, which we call the {\em associated} orientation, is given by $\ast_h$. 
\end{cor}

\begin{cor}
If  in addition we assume that $\Om^\bullet$ is a covariant calculus over a \qhsn, then the associated integral  is closed if $V\hol$, or equivalently $V\ahol$, does not contain the trivial corepresentation as a sub-comodule. 
\end{cor}
\demo
 By Corollary \ref{closurev} we know that the associated integral is closed if $V^{2n-1}$ does not contain the trivial comodule as a sub-comodule. But by the first part of the above lemma this is equivalent to $V^1$ not containing the trivial comodule as a submodule. Finally, we note that (despite not being a morphism in $\sgmm$) the $\ast$-map  brings  copies of the trivial comodule in $\Om\hol$ to  copies of the trivial comodule in $\Om\hol$, and vice versa. Hence, one need only check either $V\hol$ or $V\ahol$. 
\qed

\begin{remark}
Note that the Hermitian condition is not necessary for the existence of a $*$-orientation, one exists for any $*$-calculus containing a symplectic form. Moreover, by dropping the factor $i^{a-b}$ from the definition of $\ast_h$, it is possible to define a Hodge map for an almost  symplectic form which is not necessarily associated to a complex structure. For a discussion of such maps in the classical case see \cite{JLB, Yau}.  
\end{remark}

\begin{remark}
The use of quantum integers in the definition of the Hodge map is motivated by their appearance in the multiplicative structure of the \hk calculus in \textsection 4.4,   in the \hk calculus of the quantum Grassmannians in \cite{MM}, and in the calculus introduced for the full quantum flag manifold of $\bC_q[SU_3]$ in \cite{SM}. It is worth stressing that the Hodge parameter need not depend on a deformation parameter: indeed, the definition of the Hodge map is well-defined for \algs which are not deformations and even for commutative algebras. As more examples of \nc Hermitian structures emerge, it can be expected that a more formal definition of the Hodge map will appear (see for example the recent paper on braided Hodge maps \cite{SMH}). 
\end{remark}

\subsection{An Hermitian Structure for the Heckenberger--Kolb Calculus over $\cpn$}

In this subsection we construct a covariant Hermitian $(\Om^{(\bullet,\bullet)},\k)$ structure for the Hecken- berger--Kolb calculus over $\cpn$. In the classical case, it follows from the classification of covariant metrics on complex projective space that $\k$ is equal, up to scalar multiple, to the fundamental form of the Fubini--Study metric.

Throughout this subsection we will, by abuse of notation, denote $\Phi(L), \Phi(\vol)$, and $\Phi(\ast_q)$, by $L,\vol$, and $\ast_q$, respectively.

 \begin{lem} A left $G$-coinvariant closed  form is given by
\bas
\k : = i \sum_{k,l =1}^{n+1}  q^{2k} \del z_{kl} \wed \adel z_{lk}  =  \unit\inv \bigg(i \sum_{a=1}^{n} 1 \oby e^+_a \wed e^-_a \bigg). 
\eas
\end{lem}
\demo
The fact that $\k$ is closed follows from Lemma \ref{v3closedform} and Lemma \ref{multformweights}.  Left $G$-coinvariance of $\k$, and equality of the two given presentations, follow from 
\bas
\unit \bigg(\sum_{k,l =1}^{n+1}  q^{2k} \del z_{kl} \wed \adel z_{lk} \bigg) = & \sum_{k,l=1}^{n+1} \sum_{a,b,c,d=1}^{n+1}   q^{2k} u^k_aS(u^b_l)u^l_cS(u^d_k) \oby [\del z_{ab}] \wed [\adel z_{cd}] \\
= &    \sum_{k =1}^{n+1} \sum_{a,b ,d=1}^{n+1}  q^{2k}  u^k_aS(u^d_k) \oby [\del z_{ab}] \wed [\adel z_{bd}]\\
= & \sum_{a=2}^{n+1} 1 \oby [\del z_{a1}] \wed \big(q^{2a} [\adel z_{1a}]\big)\ =  \sum_{a=1}^{n} 1 \oby e^+_a \wed e^-_a.\\
\eas
Note that in the penultimate line we have used the fact that $[\del z_{kl}] = [\adel z_{kl}] = 0$,  for $k,l \geq 2$ or $k = l =1$ (as presented in Lemma \ref{basislem}) and in the last line we have used the elementary identity $\sum_{k=1}^{n+1} q^{2k}  u^k_aS(u^d_k) = \d_{ad} q^{2a} 1$.
\qed

\begin{lem}
It holds that 
\bas
\unit(\k^l) := i^{l \textrm{\em \, (mod } 2)} [l]_{q}! \sum_{I \in O(l)} 1 \oby e^+_I \wed e^-_I,
\eas
where $O(l)$ is the set of all  ordered subsets of  $\{1,\ldots, n\}$ with $l$ elements. 
\end{lem} 
\demo
Assuming that the proposition holds for $l$, we have
\bas
\unit(\k^{l+1})  = & \unit(\k) \wed \Big( i^{l \text{ (mod 2}) }[l]_{q}! \sum_{I \in O(l)} 1 \oby e^+_I \wed e^-_I \Big)\\ 
               = &   i^{l \text{ (mod 2}) + 1} [l]_{q}!  \sum_{I \in O(l)} \sum_{a=1}^n 1 \oby e^+_I \wed e^+_i \wed e^-_i \wed e^-_I.
\eas
In order to re-express this identity, we introduce the following notation: For $I \in O(l+1)$, denote by $I_a$, and $_aI$, the $(l+1)$-tuples where the $a\th$-entry  has been bubbled through to the last, respectively first, position. Now, as a little thought will confirm, it holds that 
\bas
\unit(\k^{l+1})   & =   i^{l \text{ (mod 2}) + 1} [l]_{q}!  \sum_{I \in O(l+1)} \sum_{a=1}^{l+1} 1 \oby e^+_{I_a} \wed e^-_{_aI}.
\eas 
The commutation relations of the calculus imply that
$
e^+_{I_a} \wed e^-_{_aI} = (-1)^lq^{l-2(a-1)} e^+_{I} \wed e^-_{I}.
$
Hence
\bas
\unit(\k^{l+1})            =   & \,  i^{l \text{ (mod 2}) + 1} (-1)^l [l]_{q}!  \sum_{I \in O(l+1)} (q^{l} + q^{l - 2} + \cdots + q^{-l + 2} + q^{-l}) 1 \oby e^+_I \wed e^-_I \\
          =    &  \, i^{l + 1\text{ (mod 2})}[l]_q! [l+1]_q \sum_{I \in O(l+1)} 1\oby e^+_I \wed e^-_I  \\
=  &\,  i^{l + 1\text{ (mod 2})} [l+1]_q! \sum_{I \in O(l+1)} 1 \oby e^+_I \wed e^-_I. 
\eas
Finally, since the proposition clearly holds for  $l=1$, we can  conclude that it holds for all $l \in \bN$.
\qed

 \begin{prop} The pair $(\Om^{(\bullet,\bullet)},\k)$ is a covariant Hermitian structure for the Heck-enberger--Kolb calculus over $\cpn$. 
\end{prop}
\demo
Since $\k$ is closed, it follows from Lemma \ref{closedeqcentral} that it is central. That $\k$ is real follows from 
\bas
\k^* : = & -i \sum_{k,l =1}^n q^{2k} \big(\del z_{kl} \wed \adel z_{lk}\big)^* = i \sum_{k,l =1}^n  q^{2k}(\adel z_{lk})^* \wed (\del z_{kl})^* \\
=  & \, \, \, i \sum_{k,l =1}^n  q^{2k} \del z_{kl} \wed \adel z_{lk} = \k.
\eas

It remains to show that $L^{n-k}: V^k \to V^{2n - k}$ is an isomorphism, for all $k=0, \ldots, n-1$. Since Lemma \ref{CalcBasis} shows that $\dim(V^k) = \dim(V^{2n-k})$, we need only show that  $L^{n-k}$ has zero kernel in $V^k$. To this end, consider the  decomposition $V^k \simeq \bigoplus_{r \geq 0} V^k_r $ where
\bal \label{repdecomp}
V^k_r := \{e^+_I \wed e^-_J \in V^k \,|\, |I \cap J|  = r \}.
\eal  
For  $v \in  V^{k} \cap \ker(L^{n-k})$, denote its decomposition \wrt  (\ref{repdecomp}) by $v = \sum_{r=1}^m v_r$, where without loss of generality we assume that $v_m \neq 0$. Since $L^{n-k+m}$ acts as zero on $V^k_r$, for $r < m$, we have   
\bas
L^{n-k + m}(v) = \sum_{r=0}^{m} L^{n-r + m}(v_r) = L^{n-r + m}(v_m). 
\eas
Hence, the proposition would follow if we could show that $L^{n-r + m}$ had trivial kernel in  $V^k_m$. But this follows from the fact that, for any $e^+_I \wed e^-_J \in V^k_r$, there exists a $m \in \bZ$, \st
\bas
L^{n-k+r}(e^+_I \wed e^-_J) = \pm q^{m}e^+_{I\cup (I \cup J)^c} \wed e^-_{J \cup (I \cup J)^c}, 
\eas
where $\cup$ denotes set union followed by reordering. 
\qed

\begin{cor}
Denoting $e^\bullet : = e^+_1 \wed \cdots e^+_n \wed e^-_1 \wed \cdots e^-_n$, it holds that
\bas
\vol(e^\bullet) =  i^{-n \text{ mod } 2}.
\eas
\end{cor}
\demo
This follows from the calculation 
\bas
~~~~~~~~~ 1 = \vol(\ast_q(1)) = \vol\Big(\frac{1}{[n]_q}\k^n\Big) = i^{n \text{ mod } 2}\vol\Big(\frac{[n]_q}{[n]_q}e^\bullet\Big) =  i^{n \text{ mod } 2} \vol(e^\bullet). ~~~~~~~~~~~ \text{\qed}
\eas

\begin{lem}\label{uniquenessHF}
Up to scalar multiple $\k$ is the unique coinvariant Hermitian form in $\Om^{(1,1)}$. 
\end{lem}
\demo  Since the calculus is factorisable we have $V^{(1,1)} \simeq V\hol \oby V\ahol$. Using an elementary representation theoretic argument it can be shown that the decomposition of $V\hol \oby V\ahol$ into irreducible summands contains a unique copy of the trivial comodule $\bC$. Thus, any other coinvariant $(1,1)$-form is a scalar multiple of $\k$. \qed

We will now look at Lefschetz decomposition and the associated Hodge map for some low dimensional examples. 

\begin{eg}
For $\bC_q[\bC P^1]$, we have 
\bas
\ast_h(e^+_1) = -ie^+_1, & &  \ast_h(e^-_1) = ie^-_1, & &  \ast_h(1) = \k = ie^+_1\wed e^-_1.
\eas 
\end{eg}

\begin{eg}\label{egLH}
For $\bC_q[\bC P^2]$, Lefschetz decomposition is trivial except for $V^{(1,1)}$, where it reduces to  $V^{(1,1)} \simeq L(1) \oplus P^{(1,1)}$. The fact   $e^+_1 \wed e^-_2, e^+_2 \wed e^-_1 \in P^{(1,1)}$ follows from
\bas
L(e^+_1 \wed e^-_2) = & e^+_1 \wed \k \wed e^-_2 =  i e^+_1 \wed (e^+_1 \wed e^-_1  + e^+_2 \wed e^-_2) \wed e^-_2  = 0,
\eas
and the analogous calculation for $L(e^+_1 \wed e^-_2)$.  Moreover, 
\bas
L(e^+_1 \wed e^-_1 - q^{-2} e^+_2 \wed e^-_2) = & \, e^+_1 \wed \k \wed e^-_1 - q^{-2} e^+_2 \wed \k \wed e^-_2\\
 = & \, i e^+_1 \wed e^+_2 \wed e^-_2 \wed e^-_1 - i q^{-2} e^+_2 \wed e^+_1 \wed e^-_1 \wed e^-_2\\
= & - iq\inv  e^+_1 \wed e^+_2 \wed e^-_1 \wed e^-_2 + i q \inv e^+_1 \wed e^+_2 \wed e^-_1 \wed e^-_2\\ 
= & \, 0
\eas
implies that  $e^+_1 \wed e^-_1 - q^{-2} e^+_2 \wed e^-_2 \in P^{(1,1)}$. Since the set
\bas
 \{\k, e^+_1 \wed e^-_2, e^+_2 \wed e^-_1, e^+_1 \wed e^-_1 - q^{-2} e^+_2 \wed e^-_2 \}
\eas
is clearly a basis for $V^{(1,1)}$, we must have that
\bas
P^{(1,1)} = \spn_{\bC}\{e^+_1 \wed e^-_2, e^+_2 \wed e^-_1, e^+_1 \wed e^-_1 - q^{-2} e^+_2 \wed e^-_2 \}.
\eas

Setting the Hodge parameter equal to $q$, the action of the associated Hodge map on $V^1$ is given by 
\bas
\ast_q(e^+_1) =   e^+_1\wed e^+_2 \wed e^-_2, ~~~~~& & \ast_q(e^+_2) = -  q e^+_1\wed e^+_2 \wed e^-_1, \\
 \ast_q(e^-_1) =   q\inv e^+_2\wed e^-_1 \wed e^-_2, & & \ast_q(e^-_2) = -  e^+_1\wed e^-_1 \wed e^-_2.~
\eas
The action of $\ast_h$ on the primitive elements of $V^2$ is  given by 
\bas
\ast_q|_{P^{(2,0)}} = \id, & & \ast_q|_{P^{(1,1)}} = - \id, & & \ast_q|_{P^{(0,2)}} = \id.
\eas
\end{eg}

\subsection{A Conjectured Hermitian Structure for $\bC_q[G/L_S]$}

Let  $\frak{g}$ be a complex semisimple Lie \alg of rank $r$ and $U_q(\frak{g})$ the corresponding Drinfeld--Jimbo quantised enveloping algebra \cite[\textsection 6.1]{HK}.  For $S$ a subset of simple roots,  denote by $\pi_S:\bC_q[G] \to \bC_q[L_S]$  the Hopf \alg map dual to the inclusion $U_q(\frak{l}_S) \hookrightarrow U_q(\frak{g})$, where
\bas
U_q(\frak{l}_S) := \la K_i, E_j, F_j \,|\, i = 1, \ldots, r; j \in S \ra.
\eas 
The quantum homogeneous space of this coaction is called the {\em quantum flag manifold} corresponding to $S$, and is denoted by $\bC_q[G/L_S]$. (See \cite{HK,HKdR} for a more detailed presentation of this definition.)

If $S = \{1, \ldots, r\}\bs \a_i$ where $\a_i$ appears in any positive root with coefficient at most one, then we say that the quantum flag manifold is {\em irreducible}. It follows that $\cpn$ is an irreducible quantum flag manifold. Moreover, Theorem \ref{HKClass} holds  for this larger family of quantum homogeneous spaces.

\begin{thm} \cite[\textsection 2]{HK}
There exist exactly two non-isomorphic irreducible left-\linebreak covariant first-order differential calculi of finite dimension over $\bC_q[G/L_S]$. 
\end{thm}

In \cite{HKdR} the maximal prolongation of the direct sum of these two calculi is shown to have a unique covariant complex structure $\Om^{(\bullet,\bullet)}$. Using a representation theoretic argument directly analogous to that in Lemma \ref{uniquenessHF}, it can be shown that the each  $\Om^{(1,1)}$ contains a left-coinvariant  form $\k$ that is unique up to scalar multiple.

\begin{conj}
For every irreducible quantum flag manifold $\bC_q[G_0/L_0]$, the pair \linebreak $(\Om^{(\bullet,\bullet)},\k)$ is a covariant Hermitian structure for the \hk calculus. 
\end{conj}

\section{Metrics, Inner Products, and Operator Adjoints}

In the previous section, Hermitian structures were introduced as a \nc generalisation of the fundamental form of an Hermitian metric, and an associated  Hodge map was defined through the classical Weil formula. In this section we bring this series of ideas full circle by defining a metric through the classical definition of the Hodge map. This allows for the construction of adjoint operators for all $G$-comodules maps, which is one of the principal motivations of the paper and an indispensable tool in \textsection 6 and \textsection 7. An interesting new phenomenon to emerge is a deformation of the classical Lefschetz identities  to a representation of the quantised enveloping \alg of $\frak{sl}_2$, see Corollary 5.14. 

Throughout this section $\Om^\bullet$ denotes a differential $*$-calculus of total dimension $2n$, and $(\Om^{(\bullet,\bullet)},\k)$ denotes an Hermitian structure for $\Om^\bullet$.

\subsection{Metrics}

By reversing the classical order of definition, we use the Hodge map to associate a metric to any Hermitian structure. 

\begin{defn}
The {\em metric} associated to the Hermitian structure $(\Om^{(\bullet,\bullet)},\k)$ is defined to be the map $g:\Om^\bullet \oby_M \Om^\bullet \to M$  for which $g(\Om^k \oby_M \Om^l) = 0$, for all $k \neq l$, and 
\bas
g(\w \oby \nu) = \vol(\w \wed \ast_h(\nu^*)), & &  \w,\nu \in \Om^k.
\eas
\end{defn}

The  $\bN^2_0$-decomposition, and the Lefschetz decomposition, of the de Rham complex of a classical Hermitian manifold are orthogonal \wrt the metric  \cite[Lemma 1.2.24]{HUY}. We now show that this carries over to the \nc setting. Moreover, we show as a consequence that the metric is conjugate symmetric.

\begin{lem}\label{orthogwrtg} It holds that
\begin{enumerate}  

\item the $\bN^2_0$-decomposition of $\Om^\bullet$ is orthogonal \wrt $g$,

\item the Lefschetz decomposition of  $\Om^{\bullet}$ is orthogonal \wrt $g$.
\end{enumerate}
\end{lem}
\demo 
\begin{enumerate}

\item The first part of Lemma \ref{Hodgeprop} implies that, given any $\w \in \Om^{(a,b)},\nu \in \Om^{(a',b')}$, for which  $a+b=a'+b'$ but $(a,b) \neq (a',b')$, then the product $\w \wed \ast_h(\nu^*) \notin \Om^{(n,n)}$. It now follows from Lemma \ref{diamond}  that $g(\w \oby_M \nu) = 0$.

\item  For $\a \in P^k, \b \in P^l$, orthogonality of the $\bN_0$-grading  implies that    
$
g\big(L^i(\a) \oby_M L^j(\b)\big)
$
is nonzero only if  $2i+k = 2j + l =: m$. Assuming that $\b \in P^{(a,b)} \sseq P^l$, we have
\bas
g\big(L^{\frac{1}{2}(m-k)}(\a) \oby_M L^{\frac{1}{2}(m-l)}(\b)\big) = & \, \vol\big(L^{\frac{1}{2}(m-k)}(\a) \wed \ast_h L^{\frac{1}{2}(m-l)}(\b^*)\big) \\
= &  \l\vol\big(L^{\frac{1}{2}(m-k)}(\a) \wed L^{n-\frac{1}{2}(m+l)}(\b^*)\big)\\
= & \l \vol\big(L^{n-\frac{1}{2}(k+l)}(\a) \wed \b^*\big),
\eas
where $\l := (-1)^{\frac{l(l+1)}{2}}i^{a-b}\frac{[\frac{1}{2}(m-l)]_h!}{[n-\frac{1}{2}(m+l)]_h!}$.
Assuming now that $k<l$, which is to say that $l = k+r$, for some $r \in 2\bN_{>0}$, 
\bas
\l\inv g\big(L^{\frac{1}{2}(m-k)}(\a) \oby_M L^{\frac{1}{2}(m-l)}(\b)\big) = & \, \vol\big(L^{n - k + \frac{r}{2}}(\a) \wed \b^*\big).
\eas
Since $\a \in P^k$, we must have $L^{n-k+\frac{r}{2}}(\a) = 0$. The proof for $k>l$ is analogous.   \qed
\eet

\begin{cor}\label{conjsym}
It holds that $g(\w \oby_M \nu) = (g(\nu \oby_M \w))^*$, for all $\w, \nu \in \Om^\bullet$.
\end{cor}
\demo
By the above lemma, it suffices to prove the result for $g(L^j(\a) \oby_M L^j(\b))$, for some $\a, \b \in P^{(a,b)} \sseq P^k$. This is done as follows
\bas
~~~~g(L^j(\a) \oby_M L^j(\b)) = & \, \vol(L^j(\a) \wed \ast_h(L^j(\b^*))) \\
= & \, (-1)^{\frac{k(k+1)}{2}}i^{b-a}\frac{[j]_h!}{[n-k-j]_h!}  \vol\big(L^j(\a) \wed  L^{n-k-j}(\b^*)\big)\\
= & \, \Big((-1)^{\frac{k(k+1)}{2}} i^{a-b}\frac{[j]_h!}{[n-k-j]_h!} (-1)^{k^2}\vol\big(L^j(\b) \wed  L^{n-k-j}(\a^*)\big)\Big)^*\\
= & \Big((-1)^{\frac{k(k+1)}{2}}i^{b-a}\frac{[j]_h!}{[n-k-j]_h!} \vol\big((L^j(\b) \wed L^{n-k-j}(\a^*)\big)\Big)^*\\
~~~~ = & \,\big(g(L^j(\b) \oby_M L^j(\a))\big)^*. ~~~~~~~~~~~~~~~~~~~~~~~~~~~~~~~~~~~~~~~~~~~~~~~~~ \text{\qed}
\eas

\subsection{Inner Products and  Operator Adjoints}

In this subsection we specialise to the case where $\Om^{\bullet}$ is a covariant calculus over a \qhs $M$, and $(\Om^{(\bullet,\bullet)},\k)$ is a covariant Hermitian structure. Following the classical order of definition, we introduce an inner product by composing the associated metric with the Haar functional. In order for this to well-defined, however, we need to impose a positive definiteness condition on our Hermitian structure.

\begin{defn}
An Hermitian structure is said to be {\em positive definite} if an inner product is given by
\bal \label{PosDef}
\la\cdot,\cdot\ra_V:  V^{\oby 2} \to \bC, & & [\w] \oby [\nu] \mto [g(\w \oby_M \nu)].
\eal  
\end{defn}

Note that $\la\cdot,\cdot\ra$ is well-defined because of our assumption that $\Om^\bullet \in \sgmm$. Another point, which is easily checked, is that positive definiteness of an Hermitian structure is independent of the choice of Hodge parameter. 

In general it can be quite tedious to verify positive definiteness; the following lemma shows us that we need only do so for primitive elements.
\begin{lem}\label{Gaussbin}
For $\a,\b \in P^\bullet$, and $\binom{a}{b}_h := \frac{[a]_h!}{[b]_h![a-b]_h!}$  the Gaussian binomial coefficient, it holds that  
\bas
\la L^j(\a),L^j(\b)\ra_V = \binom{n-j-k}{j}_h\inv \la \a,\b \ra.
\eas
\end{lem}
\demo
Assuming, without loss of generality, that $\b \in P^{(a,b)} \sseq P^k$, we have 
\bas
~~~~~~ \la L^j(\a),L^j(\b)\ra_V = &\, \vol\big(L^j(\a) \wed \ast_h(L^j(\b^*))\big)\\
 = & (-1)^{\frac{k(k+1)}{2}}i^{b-a} \frac{[j]_h!}{[n-j-k]_h!} \vol\big(L^j(\a)\wed L^{n-j-k}(\b^*) \big)\\
= & \frac{[j]_h![n-k]_h!}{[n-j-k]_h!}\vol(\a\wed \ast_h(\b^*))
=  \binom{n-j-k}{j}_h\inv \la\a,\b\ra_V. ~~~~~~~~~ \text{\qed}
\eas

\begin{cor}
If $\la\cdot,\cdot\ra_V$ restricts to an inner product on the space of primitive elements, then it is an inner product on all of $V^\bullet$.
\end{cor}

We are now ready to introduce the inner product associated to an Hermitian structure and to establish the existence of adjoints \wrt  this pairing.

\begin{lem}
For $\ast_h$ the Hodge map of a positive definite Hermitian structure, an inner product is given by
\bal \label{fiberinnerprod}
\la \cdot,\cdot\ra: \Om^\bullet \oby \Om^\bullet \to \bC, & & \w \oby \nu \mto   \int \w \wed \ast_h(\nu^*) = \haar(g(\w \oby_M \nu^*)).
\eal
Moreover, the Peter--Weyl decomposition of $\Om^\bullet$ is orthogonal \wrt $\la \cdot,\cdot\ra$.
\end{lem}
\demo
By Corollary \ref{conjsym}, we need only establish positive definiteness. Let $\big\{[\w_k]\big\}_k$, for $\w_k \in \Om^{\bullet}$, be an \onb of $V^\bullet$ \wrt the inner product (\ref{PosDef}). In what follows we denote  $\unit(\w) =:  \sum_{} f_{k} \oby [\w_k]$, and tacitly assume the isomorphism  
$
\id \oby \e = \unit \inv:G \coby \Phi(M) \to M.
$
Noting that a morphism is given by $\ol{g} := (\id \oby *)g$, positive definiteness of the bilinear form follows from 
\bas
\big<\w, \w^*\big> = &\, \haar \circ (\id \oby \Phi(\ol{g})) \circ \unit \big(\w \oby_M \w^*\big) = \sum_{k,l} \haar(f_k f_l^*)  [g(\w_k \oby_M \w_l)]  \\
= & \,  \sum_{k}  \haar(f_k f_k^*) \in \bR_{>0}.
\eas
Orthogonality of the Peter--Weyl decomposition of $\Om^\bullet$ is established similarly. \qed


\begin{cor}\label{EqAdj}
Any left $G$-comodule map $f:\Om^\bullet \to \Om^\bullet$ is adjointable \wrt $\la\cdot,\cdot\ra$.  Moreover, if $f$ is self-adjoint, then it is diagonalisable, and commuting diagonalisable maps are simultaneously diagonalisable. {\bf }
\end{cor}
\demo
Since $f$ is a left $G$-comodule map  $f(\Om^\bullet_V) \sseq \Om^\bullet_V$, for all $V \in \wh{G}$. Adjointability of $f$ now follows from finite-dimensionality of $\Om^\bullet_V$ and the fact that Peter--Weyl decomposition is orthogonal \wrt  $\la\cdot,\cdot\ra$. Analogously,  $f$ can be shown to be diagonalisable whenever it is self-adjoint, and so, commuting diagonalisable maps can be shown to be simultaneously diagonalisable.
\qed

\begin{remark}
In \cite[\textsection 3]{KMT} a calculus $\Om^\bullet$ is defined to be {\em non-degenerate}  if, whenever $\w \in \Om^k(M)$, and $\w \wed \nu = 0$, for all $\nu \in \Om^{n-k}(M)$, then necessarily $\w = 0$. Clearly, the existence of a positive definite Hermitian form for a $*$-calculus implies non-degeneracy.  
\end{remark}

\subsection{Examples of Operator Adjoints}

We now consider three explicit examples of adjointable operators: the Hodge map, the Lefschetz operator, and the differentials $\exd,\del,\adel$. The Hodge operator is shown to be unitary, while the adjoints of the other operators are shown to admit explicit descriptions in terms of the Hodge map. In the case of the Lefschetz map, this allows us to establish a $h$-deformation of the classical Lefschetz identities. 

Note that throughout this subsection,  we continue to assume that $\Om^{\bullet}$ is a covariant $*$-calculus over a \qhs $M$, and $(\Om^{(\bullet,\bullet)},\k)$ is a covariant Hermitian structure. Moreover, $(\Om^{(\bullet,\bullet)},\k)$ is assumed to be positive definite. To avoid confusion with the $*$-map, the symbol $\dagger$ will be use to denote the adjoint of an operator.

\subsubsection{Unitarity of the Hodge Map}

Here we show that, just as in the classical case,  $\ast_h$ is unitary. (Note that this property is assumed in the definition of the \nc Hodge map in \cite[Definition 5.20]{FGR}.)

\begin{lem}
For all values of the Hodge parameter $h$, the Hodge map is unitary.
\end{lem}
\demo
For $\a,\b \in P^{(a,b)} \sseq P^k$, and $j\geq 0$, we have 
\bas
 \la\ast_h(L^j(\a)),\ast_h(L^j(\b))\ra &  =  \int \ast_h(L^{j}(\a)) \wed \ast_h^2(L^j(\b^*)) \\
& = (-1)^{\big(\frac{k(k+1)}{2} + k\big)}i^{a-b}\frac{[j]_h!}{[n-j-k]_h!} \int L^{n-j-k}(\a) \wed L^j(\b^*) \\
& = (-1)^{\frac{k(k+1)}{2}} i^{b-a}\frac{[j]_h!}{[n-j-k]_h!} \int L^{j}(\a) \wed L^{n-j-k}(\b^*)\\
& = \int L^j(\a) \wed \ast_h(L^j(\b^*))  = \la L^j(\a), L^j(\b)\ra. 
\eas
The result now follows from orthogonality of the Lefschetz decomposition. \qed

\subsubsection{The Dual Lefschetz Operator and the Lefschetz Identities}

We now present an explicit formula for the adjoint of $L$ in terms of $\ast_h$, this is again a direct generalisation of a well-known classical formula \cite[Lemma 1.2.2 3]{HUY}.

\begin{lem}\label{adjointability}
It holds that $\Lambda := L^\dagger = \ast_h \inv  L \, \ast_h$.
\end{lem}
\demo
For $\w,\nu \in \Om^k$, we have
\bas
~~~~~~~~~~~~~~~~~~ \big<L(\w),\nu \big> & = \int L(\w) \wed \ast_h(\nu^*) = \int \w  \wed L \ast_{h}(\nu^*) \\
 & = \int \w \wed \ast_h  \big(\ast_h \inv L \ast_h(\nu^*)\big)   = \la \w,\ast_h \inv L \ast_{h}(\nu)\ra. ~~~~~~~~~~~~~~~~~~ \text{ \qed}
\eas

Classically the primitive forms are defined to be those contained in the kernel of $\Lambda$. The following corollary derives this as a consequence of our definition of primitive forms.

\begin{cor}
It holds that $P^k = \ker(\Lambda:\Om^k \to \Om^{k-2})$
\end{cor}
\demo
For $\a \in P^{(a,b)} \sseq P^k$, the inclusion $P^k \sseq \ker(\Lambda:\Om^k \to \Om^{k-2})$ follows from
\bas
\Lambda(\a) = & \ast_h\inv L \ast_h(\a) = (-1)^{\frac{k(k+1)}{2}}i^{a-b}\frac{1}{[n-k]_h!}\ast_h L^{n-k+1}(\a) = 0.
\eas
For the opposite inclusion consider, for $j>0$,
\bas
0 = \Lambda(L^j(\a)) = \ast_h \inv L \ast_h (L^j(\a))= (-1)^{\frac{k(k+1)}{2}}i^{a-b}\frac{[j]_h!}{[n-j-k]_h!}\ast_h\inv L^{n-j-k+1}(\a).
\eas
Since $\ast_h$ is an isomorphism, we must have $L^{n-j-k+1}(\a) = 0$, and so, that  $\a \in P^{k+j}$. Since $P^k \cap P^{k+j} = 0$,  we must $\a = 0$.
\qed

\bigskip

Consider now  the {\em counting operators} 
\bas
H, K : \Om^\bullet \to \Om^\bullet, & & H(\w) = (k-n) \w, & & K(\w) = h^{k-n} \w, & &   \w \in \Om^k.
\eas
For a classical Hermitian manifold the operators $H$, $L$, and $\Lambda$, define a representation of $\frak{sl}_2$ \cite[Proposition 1.2.26]{HUY}. We now show that in the \nc setting $H,L,\Lambda$, and $K$ give a representation of the quantised enveloping \alg of $\frak{sl}_2$.

\begin{prop}\label{LIDS}
We have the relations
\bas
[H,L]_{h^{-2}} = [2]_h LK, & & [L,\Lambda] = H, & & [H,\Lambda]_{h^2} = - [2]_{h^2} K \Lambda,
\eas
where $[A,B]_{h^{\pm 2}} = A B - h^{\pm 2}B A$.
\end{prop}
\demo
Beginning with the first relation, for $\w \in \Om^k$,
\bas
[H,L]_{h^{-2}}(\w) = & H L(\w) - h^{-2}LH(\w) =  \big([k+2-n]_h - h^{-2}[k-n]_h\big)L(\w)\\
                                  = & \big(h^{k-n}[2]_h + h^{-2}[k-n]_h  - h^{-2}[k-n]_h\big)L(\w)\\
                                  = &\, h^{k-n}[2]_h L(\w) = [2]_h L K(\w).
\eas 
Noting that $H$ and $K$  are self-adjoint operators, we see that the third relation is  the operator adjoint of the first.

Coming finally to the second relation, for $\a \in P^{(a,b)} \sseq P^k$, we have
\bas
L\Lambda(L^j(\a)) =    L \ast_{h}\inv L \ast_{h}(L^j(\a))  =  & \, L  \ast_{h} \inv  L\big((-1)^{\frac{k(k+1)}{2}}i^{a-b}\frac{[j]_h!}{[n-j-k]_h!} L^{n-j-k}(\a) \big)\\
 =    &  (-1)^{\frac{k(k+1)}{2}}i^{a-b}\frac{[j]_h!}{[n-j-k]_h!}L \ast_{h}\inv L^{n-j-k+1}(\a)\\ 
 =    &  (-1)^{\frac{k(k+1)}{2}+k}i^{a-b}\frac{[j]_h!}{[n-j-k]_h!} L  \ast_h L^{n-j-k+1}(\a)\\
 =    & [j]_h[n-j-k+1]_h L^j(\a). 
\eas  
Similarly, it can be shown that 
\bas
\Lambda L(L^j(\a)) = [j+1]_h  [n-j-k]_h L^j(\a).
\eas
Hence,
\bas
~~~~~~~~~~~~~  [L,\Lambda]L^j(\a)
 = &  \big([j]_h[n-j-k+1]_h  - [j+1]_h [n-j-k]_h\big) L^j(\a)\\
                                    = &  \big(h\inv [j]_h [n-j-k]_h  + h^{n-j-k}[j]_h \\
 &~~~ ~~ - h\inv[j]_h[n-j-k]_h - h^j[n-j-k]_h\big)L^j(\a)\\
= & \, [2j+k-n]_h L^j(\a)  = H L^j(\a). ~~~~~~~~~~~~~~~~~~~~~~~~~~~~~~~~~~~~~~~~~\text{\qed}
\eas

Clearly, for $h = 1$, we get a representation of the Lie \alg $\frak{sl}_2$. For the case of $h \neq 1$, we get a representation of the  quantised universal enveloping \alg  of $\frak{sl}_2$ (where we use the conventions presented in \cite[\textsection{3.1.1}]{KSLeabh}).

\begin{cor}
A representation $\rho$ of $U_h(\frak{sl_2})$   is given by 
\bas
\rho(E) = L, & & \rho(K) = K, & & \rho(F) = \Lambda.
\eas
\end{cor}
\demo
It is clear that $\rho(KK\inv) = \rho(K\inv K) =1$. Moreover, 
\bas
\rho(KEK\inv)  = & KLK\inv = h^2L = h^2\rho(E), 
\eas
and
\bas
 \rho(KFK\inv) = & K\Lambda K\inv = h^{-2}L = h^{-2}\rho(E).
\eas
Finally, for $\w \in \Om^k$, we have
\bas
~~~~~~~~~~~~~~~~~ \rho([E,F])(\w) = & [L,\Lambda](\w) = H(\w) = [k-n]_{h}(\w) =  \frac{h^{k-n}-h^{-(k-n)}}{h-h\inv}\w \\ = &  \frac{K-K\inv}{h - h \inv}(\w) = \rho\Big(\frac{K-K\inv}{h - h \inv}\Big)(\w). ~~~~~~~~~~~~~~~~~~~~~~~~~~~~~~~\text{\qed}
\eas

Finally, we describe the irreducible representations of  $U_h(\frak{sl}_2)$. Note that by taking appropriate unions, the Lefschetz decomposition can be reproduced from this decomposition. In fact,  this is how the Lefschetz decomposition is established classically.

\begin{lem}
The irreducible representations of $U_h(\frak{sl}_2)$  are given by
\bas
\bigoplus_{j \geq 0} L^j(\a),  & & \a \in P^{(a,b)}, \, \, (a,b) \in \bN^2_0.
\eas
\end{lem}

\subsubsection{Codifferential Operators}

 We call the adjoints of $\exd,\del$, and $\adel$   the {\em codifferential}, {\em holomorphic codifferential}, and {\em anti-holomorphic codifferential}, respectively.  Classically, these operators have expressions in terms of the Hodge operator analogous to the expression given above for the dual Lefschetz operator. The following lemma shows that this is also true in the \nc setting.

\begin{lem}
It holds that
\bas
\exd^\dagger = - \ast_h  \exd  \ast_h, &&  \del^\dagger =   - \ast_h  \adel \ast_h,&& \adel^\dagger =   - \ast_h  \del  \ast_h.
\eas
\end{lem}
\demo
For $\w \in \Om^k, \nu \in \Om^\bullet$,  the Leibniz rule  and  closure of the integral imply that
$
{0 = \int \adel(\w \wed \nu)}$ ${= \int \adel \w \wed \nu + (-1)^k \int \w \wed \adel \nu},
$
and so,
$
\int \adel \w \wed \nu  = (-1)^{k+1} \int \w \wed \adel\nu.
$
This in turn implies that
\bas
 \big< \w, \ast_h \, \del \ast_h(\nu)\big>  = & \int \w \wed \ast_h\big((\ast_h \del \ast_h (\nu))^*\big) =  \int \w \wed \ast_h\big(\ast_h \adel \ast_h (\nu^*)\big)\\
  = & \, (-1)^{k} \int \w \wed  \adel \ast_h (\nu^*) =  - \int \adel \w \wed  \ast_h (\nu^*)\\
                                                = & - \big<\adel \w,\nu\big>.
\eas  
Hence,  $\adel^\dagger =   - \ast_h  \del \, \ast_h$. The identities for $\exd^\dagger$ and $\adel^\dagger$ are established similarly. 
\qed

\begin{cor}\label{codstarcomm}
For all $\w \in \Om^\bullet$, it holds that
\bal \label{hodgecod}
\exd^\dagger(\w^*) = \big(\exd^\dagger(\w)\big)^*, & & \del^\dagger(\w^*) = \big(\adel^\dagger(\w)\big)^*, & & \adel^\dagger(\w^*) = \big(\del^\dagger(\w)\big)^*.
\eal
\end{cor}
\demo
This follows  from the given formulae for the codifferentials, the fact that the $*$-map commutes with the Hodge map, and identities given in (\ref{stardel}).
\qed   

\subsection{Positive Definiteness for the Heckenberger--Kolb Calculus}

We begin by directly verifying positive definiteness  of $\k$ in the two simplest cases $\bC_q[\bC P^1]$  and $\bC_q[\bC P^2]$. Throughout, by abuse of notation, we will write $\ast_q$ for $\Phi(\ast_q)$.

\begin{eg}
In this example we will verify positiveness for $\bC_q[\bC P^1]$. By Lemma \ref{Gaussbin} we only need to show positiveness on non-trivial primitive elements, which by definition are all contained in $V^1$. For $V\ahol$ it holds that 
\bas
\la e^-_1, e^-_1\ra_V = q^4\vol(e^-_1 \wed \ast_q(e^+_1)) = - iq^4\vol(e^-_1 \wed e^+_1) = iq^6\vol(e^+_1 \wed e^-_1) = q^6.
\eas
Similarly, it can be shown that $\la e^+_1, e^+_1\ra_V = q^4$. Orthogonality of $e^+_1$ and $e^-_1$ follows from 
\bas
\la e_1^+, e^-_1\ra_V = q^4 \vol(e^+_1 \wed \ast_q(e^+_1)) = - iq^4\vol(e^+_1 \wed e^+_1) = 0, 
\eas
and the analogous calculation for $\la e_1^-,  e^+_1\ra_V$. Hence $\la\cdot,\cdot\ra_V$ is indeed positive definite.
\end{eg}

\begin{eg}
We now turn to $\bC_q[\bC P^2]$. By Lemma \ref{Gaussbin} we only need to show positive definiteness on non-trivial primitive elements. By definition these elements are all contained in $V^1$ and $V^2$. For $V^1$, we have 
\bas
\|e^+_1\|_V := \la e^+_1,e^+_1\ra_V  = & \, \vol(e^+_1 \wed \ast_q((e^+_1)^*)) = q^{-4}\vol(e^+_1 \wed \ast_q(e^-_1)) 
\\
=  & \, q^{-5}\vol(e^+_1 \wed e^+_2 \wed e^-_1 \wed e^-_2 ) = q^{-5},
\eas
and similarly 
$
 \|e^+_2\|_V = q^{-5},  \|e^-_1\|_V = q^{7},$ and   $\|e^-_2\|_V = q^{9}$. 
Orthogonality of the spaces $V\hol$ and $V\ahol$ follows from Lemma \ref{orthogwrtg}. For $e^+_1, e^+_2$, we have
\bas
\la e^+_1,e^+_2 \ra_V  = \vol(e^+_1 \wed \ast_q(e^-_2)) = -i \vol(e^+_1 \wed \k \wed e^-_2)  = 0,
\eas 
and similarly that $\la e^+_2,e^+_1 \ra_V = \la e^-_1,e^-_2 \ra_V = \la e^-_2,e^-_1 \ra_V= 0$.

For $P^{(2,0)} = V^{(2,0)}$ and $P^{(0,2)} = V^{(0,2)}$ we have 
$
\|e^+_1 \wed e^+_2\|_V = q^{-11}$ and $\| e^-_1 \wed e^-_2\|_V = q^{17}$ 

Finally, for $P^{(1,1)}$  the basis elements $\{e^+_1 \wed e^-_2, e^+_2 \wed e^-_1, e^+_1 \wed e^-_1 - q^{-2} e^+_2 \wed e^-_2 \}$ are easily seen to be orthogonal. Moreover, we have  
\bas
 \la e^+_1 \wed e^-_2, e^+_1\wed e^-_2\ra_V = & \, \vol(e^+_1 \wed e^-_2 \wed \ast_q((e^+_1\wed e^-_2)^*)) =  - q^{2}\vol(e^+_1 \wed e^-_2 \wed \ast_q(e^+_2\wed e^-_1))\\
= & \, q^{2}\vol(e^+_1 \wed e^-_2 \wed e^+_2\wed e^-_1) = q^3 \,\vol(e^+_1 \wed e^+_2 \wed e^-_1 \wed e^+_2)\\
= & \, q^3,
\eas
and  $\|e^+_2 \wed e^-_1\|_V = q$ and $\|e^+_1 \wed e^-_1 - q^{-2} e^+_2 \wed e^-_2\|_V = [2]_q.$  Hence $\la\cdot, \cdot\ra_V$ does indeed give us an inner product.
\end{eg}

From these examples it is easy to see that orthogonality of the basis elements $e_I$ \wrt $\la\cdot,\cdot\ra$ extends to the general $\cpn$ case, as presented in the following lemma.

\begin{lem}
For $\cpn$ the basis $\{e_I\}_I$  is orthogonal \wrt $\la\cdot,\cdot\ra_V$.
\end{lem}

Directly extending positive definiteness to the general $\cpn$ case proves more challenging and we postpone the technical details to a subsequent work. However, using a general argument, we can prove positive definiteness for the case of $q$ contained in a certain open interval in $\bR$ around $1$. 

\begin{lem}
There exists an open real interval around $1$, such that when  $q$ is contained in this interval, the Hermitian structure $(\Om^{(\bullet, \bullet)},\k)$ is positive definite.
\end{lem}
\demo
For   $q=1$,  $V^\bullet$ is just the exterior \alg of $V$. Hence,  we can extend the restriction  $\la\cdot,\cdot\ra_V: V^1 \oby V^1 \to \bC$, to a positive definite bilinear pairing on $V^\bullet$ using the standard determinant formula. Moreover, by Weil's formula it must coincide with $\la\cdot,\cdot\ra_V:V^\bullet \oby V^\bullet \to \bC$, which must then be positive definite.

For a general $q$, and some basis element $e_I$, consider the polynomial function in $q$ given by
\bas
g_I: \bR_{>0} \to \bC, & & q \mto \la e_I, e_I\ra_V.
\eas
As is easily checked, $g_I$ has real coefficients, and so, by an elementary continuity argument, there exists an open interval in $\bR$ around $1$ on which it takes positive real values.  Taking the (finite) intersection of these intervals of over all basis elements gives the required interval.
\qed

\section{Hodge Theory}

In this section, Hodge decomposition \wrt $\exd, \del,$ and $\adel$, is established and shown to imply an isomorphism between cohomology classes and harmonic forms, just as in the classical case. A \nc generalisation of Serre duality is also proved. These results give us some powerful tools with which to approach questions about cohomology. Most of the material in subsections 6.1 and 6.2 are generalisations to the quantum homogeneous space setting of results proven for quantum groups in \cite{KMT}.

Throughout this section,  $\Om^{\bullet}$ denotes a covariant $*$-calculus,  of total dimension $2n$, over a \qhs $M$. Moreover, $(\Om^{(\bullet,\bullet)},\k)$ denotes  a positive definite covariant Hermitian structure \st the associated integral is closed.

\subsection{Laplacians and Harmonic Forms}

Directly generalising the classical situation, we define the $\exd$-, $\del$-, and $\adel$-Laplacians to be, respectively, 
\bas
& \DEL_{\exd} := (\exd + \exd^\dagger)^2, &  \, &\DEL_{\del} := (\del + \del^\dagger)^2,&  &\DEL_{\adel} := (\adel + \adel^\dagger)^2.&
\eas
Moreover, we define the space of {\em $\exd$-harmonic}, {\em $\del$-harmonic}, and {\em $\adel$-harmonic} forms to be, respectively,
\bas
& \H_{\exd} :=\ker(\DEL_{\exd}),& &\H_{\del} :=\ker(\DEL_{\del}),&  &\H_{\adel} :=\ker(\DEL_{\adel}).&
\eas
When $\Om^\bullet$ is a covariant calculus over a quantum homogeneous space,  $\DEL_{\exd},\DEL_{\del}$, and $\DEL_{\adel}$, are left $G$-comodule maps, and so, each space of harmonic forms is a left $G$-comodule.

With respect to the $\bN_0$-grading on the calculus, $\DEL_\exd$ is a homogeneous map of degree $0$, implying the decomposition $\H_\exd =: \bigoplus_{k \in \bN_0} \H^k_\exd$. Moreover,  $\DEL_\del$ and $\DEL_{\adel}$ are homogeneous  maps of degree $0$ \wrt the $\bN^2_0$-grading, implying the decompositions $\H_{\del} =: \bigoplus_{(a,b) \in \bN^2_0} \H^{(a,b)}_{\del}$ and  $\H_{\adel} =: \bigoplus_{(a,b) \in \bN^2_0} \H^{(a,b)}_{\adel}$. Note that $\DEL_\exd$ is not necessarily homogeneous  \wrt the $\bN^2_0$-grading, and so, such a  decomposition is not guaranteed to exist (see Corollary \ref{LaplacianEq}).

\subsection{The Hodge Decomposition}

We now come to Hodge decomposition, the principal result of this section, which allows us to prove statements about the cohomology ring $H^\bullet$ which are independent of any choice of Hermitian structure. The fact that we can prove such statements is one of the principal justifications we provide for introducing  Hermitian structures. 

\begin{lem} \label{harmcap}
It holds that 
\begin{enumerate} 
\item $\H_{\exd} \simeq \ker(\exd) \cap \ker(\exd^\dagger)$,
\item $\H_{\del} \simeq \ker(\del) \cap \ker(\del^\dagger)$, 
\item $\H_{\adel} \simeq \ker(\adel) \cap \ker(\ol{\del}^\dagger)$.
\end{enumerate} 
\end{lem}
\demo
Since $\exd + \exd^\dagger$ and $\DEL_\exd$ are commuting self-adjoint $G$-comodule maps, it follows from Corollary \ref{EqAdj} that they are simultaneously diagonalisable, and in particular that their kernels coincide.  Now since the codomains of $\exd$ and $\exd^\dagger$ are orthogonal, we must have
$
\ker(\exd + \exd^\dagger) = \ker(\exd) \cap \ker(\exd^\dagger),
$
which proves that first identity. The proofs of the other two identities are  analogous.
\qed

\begin{thm}
The following decompositions are orthogonal \wrt  $\la\cdot,\cdot\ra$
\begin{enumerate} 
\item $\Om^{\bullet}  \simeq {\cal H}_{\exd} \oplus \exd \Om^{\bullet} \oplus \exd^\dagger\Om^{\bullet}$, 
\item $\Om^{\bullet}  \simeq {\cal H}_{\del} \oplus \del\Om^{\bullet} \oplus \del^\dagger\Om^{\bullet}  $,
\item $\Om^{\bullet}  \simeq {\cal H}_{\ol{\del}} \oplus \ol{\del}\Om^{\bullet} \oplus \ol{\del}^\dagger\Om^{\bullet}$.
\end{enumerate}
\end{thm}
\demo
Since both  $\exd \exd^\dagger$ and $\exd^\dagger \exd$ are  self-adjoint left $G$-comodule maps, Corollary \ref{EqAdj} implies that they are diagonalisable. Moreover, since  $(\exd \exd^\dagger)( \exd^\dagger \exd) = 0 = (\exd^\dagger \exd)(\exd \exd^\dagger)$, they commute and  hence are simultaneously diagonalisable.

Denoting the simultaneous eigenbasis by $\{b_i\}_{i \in I}$,  let $\l_i$ and $\mu_i$ be the  eigenvalues determined by $\exd \exd^\dagger b_i = \l_i b_i$ and  $\exd^\dagger \exd b_i = \mu_i b_i$. Since  $(\exd \exd^\dagger)(\exd^\dagger \exd)=0$, we have $\l_i\mu_i = 0$, for all $i \in I$. If $\l_i \neq 0$, then $b_i = \exd (\exd^\dagger (\l_i \inv b_i)) \in \exd\Om^\bullet$. Similarly, if $\mu_i \neq 0$, then $b_i = \exd^\dagger (\exd(\mu_i \inv b_i)) \in \exd^\dagger \Om^\bullet$. Finally, if $\l_i = \mu_i = 0$, then $b_i \in \H^\bullet_\exd$. This implies that 
$
\Om^{\bullet}  = {\cal H}_{\exd} + \exd( \Om^{\bullet}) + \exd^\dagger( \Om^{\bullet} ).
$ 

We now show that this is an orthogonal decomposition. Since $\la \exd \w, \exd^\dagger \nu \ra = \la \exd^2 \w , \nu\ra = 0$, the spaces $\exd \Om$ and $\exd^\dagger \Om$ are orthogonal. Orthogonality of ${\cal H}_{\exd}$ and $ \exd \Om^{\bullet} \oplus \exd^\dagger \Om^{\bullet}$ follows from
\bas
\la \exd \w + \exd^\dagger \nu, \rho \ra = \la  \w, \exd^\dagger \rho \ra + \la \nu, \exd \rho \ra = 0, & & \w,\nu \in \Om^\bullet, \rho \in {\cal H}_{\exd}.
\eas
The other two isomorphisms are established analogously.
\qed

\bigskip

\begin{cor} \label{harmonictoclass}  
It holds that 
\bas
\ker(\exd) \simeq \H_{\exd} \oplus \exd \Om^\bullet, & & \ker(\del)  \simeq \H_{\del} \oplus \exd \Om^\bullet, & & \ker(\adel) \simeq \H_{\adel} \oplus \exd \Om^\bullet,
\eas
and so, we have the isomorphisms
\bas
\H^k_{\exd}  \to  \,H^k_{\exd}, & &  \H^{(a,b)}_{\del} \to H^{(a,b)}_{\del}, & & \H^{(a,b)}_{\adel} \to H^{(a,b)}_{\adel}.
\eas
\end{cor}
\demo
For any $\w \in \Om^\bullet$ \st $\exd \exd^\dagger \w = 0$, we have $0 = \la \exd \exd^\dagger \w, \w\ra = \la \exd^\dagger \w, \exd^\dagger \w\ra$, and so, by positive definiteness $\exd^\dagger \w = 0$, implying that $\ker(\exd) \cap \exd^\dagger \Om^\bullet = 0$. Hodge decomposition and Lemma \ref{harmcap} now imply that $\ker(\exd) \simeq \H_\exd \oplus \exd \Om^\bullet$. The other two isomorphisms are established analogously.
\qed

\begin{cor}\label{lapcommiso}
Any linear map $A:\Om^\bullet \to \Om^\bullet$ which commutes with the Laplacian $\DEL_\exd$ induces a unique map on $H^\bullet$ for which the following diagram is commutative: 
\bas
\xymatrix{ 
\H^\bullet_{\exd}   \ar[d]_{A}                        & & &   \ar[lll]_{\simeq}    H^{\bullet}_{\exd} \ar[d]^A  \\
\H^\bullet_{\exd}   \,\,   \ar[rrr]_{\simeq}                     & & & H^{\bullet}_{\exd}.\\
}
\eas  
Moreover, if $A$ restricts to an isomorphism $\H^k_{\exd} \to \H^l_{\exd}$, for some $k,l \in \bN_0$, then the corresponding map $H^k_{\exd} \to H^l_{\exd}$ is also an isomorphism. The analogous results hold for $\DEL_\del$ and $\DEL_{\adel}$. 
\end{cor}
\demo
If $A$ commutes with the Laplacian, then clearly it maps harmonic forms to harmonic forms, and so, by Hodge decomposition it induces a map on $H^\bullet$.   
Since the map $A\inv:\Om^l \to \Om^k$ must also commute with the Laplacian, we must have an inverse $A\inv:H^l_{\exd} \to H^k_{\exd}$. The proofs for $\DEL_\del$ and $\DEL_{\adel}$ are analogous.
\qed

Using this corollary, we show that the Hodge map and the $*$-map induce isomorphisms  on the cohomology ring of $\Om^\bullet$, and present some easy but interesting consequences. 

\begin{lem} \label{kcoho}
The Hodge map $\ast_h$, and the $*$-map, commute with  the Laplacian $\DEL_\exd$, and so, induce isomorphisms on $H^\bullet_{\exd}$.
\end{lem}
\demo
The fact that the $*$-map commutes with the $\DEL_{\exd}$ follows directly from Corollary \ref{codstarcomm}. For the Hodge map, note that, for any $\w \in \Om^k$, 
\bas
~~~ ~~~[\ast_h,\DEL_\exd](\w)  = & \ast_h (\exd\exd^\dagger + \exd^\dagger \exd)(\w) - (\exd\exd^\dagger + \exd^\dagger \exd)\ast_h(\w)\\
= & - \ast_h \, \exd \ast_h \exd \ast_h(\w)  -  \ast_h^2 \, \exd \ast_h \exd(\w)  +  \exd \ast_h \exd \ast_h^2(\w)  + \ast_h \, \exd \ast_h \exd \ast_h(\w)\\
= & (-1)^{2n-k} \exd \ast_h \exd(\w)  -  (-1)^k \exd \ast_h \exd(\w) = 0. ~~~~~~~~~~~~~~~~~~~~~~~~~~~~~~~~~~~ \text{\qed}
\eas

\begin{cor} 
If the cohomology ring $H^\bullet_{\exd}$ has finite dimension, then 
\bas
\dim(H^{2k+1}_{\exd}) \in 2\bN_0,  & & \text{ for all  } k = 0, \ldots, n-1.
\eas
\end{cor}
\demo
Since the $*$-map induces an isomorphism between $H^{(a,b)}_{\exd}$ and $H^{(b,a)}_{\exd}$, it implies that they have equal dimension. Evenness of $\dim(H^{2k+1}_{\exd})$ now follows from
\bas
 ~~~~~~~~~~ \dim(H^{2k+1}_{\exd}) = \sum_{i=0}^{2k + 1} \dim(H^{2k+1-i,i)}_{\exd}) = 2 \sum_{i=0}^{k} \dim(H^{(2k+1-i,i)}_{\exd}) \in 2\bN_0. ~~~~~~~~~~~~~ \text{\qed}
\eas

\begin{cor} \label{h2ncor}
It holds that  $H^{2n}_{\exd} \neq 0$.
\end{cor}
\demo
Since $\DEL_\exd(1) = 0$, we have $H^0_{\exd} \neq 0$. The result now follows from the isomorphism $\ast_h: H^0_\exd \to H^{2n}_{\exd}$.
\qed

\subsection{Serre Duality}

We finish the section with a proof of Serre duality for Dolbeault cohomology, following the standard proof in \cite[\textsection 3.2]{HUY}.  (See also \cite{MBBJS} for a discussion of Serre duality from a noncommutative algebraic geometry point of view.)

\begin{prop}
Non-degenerate pairings are given by
\bas
H^{(a,b)}_{\adel} \by H^{(n-a,n-b)}_{\adel} \to \bC, & & ([\a],[\b]) \to \int \a \wed \b,
\eas
and the analogous pairing for $\H^\bullet_{\del}$.
\end{prop}
\demo
Recalling that $\a$ and $\b$ are $\del$-closed forms and that $\int$ is assumed to be closed,  the fact the pairing is well-defined follows from  
\bas
\int (\a + \adel \w) \wed (\b + \adel \nu)  & = \int \a \wed \b + \int \adel \w \wed \b + \int \a \wed \adel \nu + \int \adel \w \wed \adel \nu\\ 
              & = \int \a \wed \b + \int \adel (\w \wed \b) + (-1)^{a+b} \int \adel (\a \wed  \nu) + \int \adel (\w \wed \adel \nu) \\
              & = \int \a \wed \b.
\eas
Next note that for any nonzero $\w \in \H^{(a,b)}_{\adel}$, the form $\ast_h(\w^*)$ is an element of $\H^{(n-a,n-b)}_{\adel}$.  Since
$
\int \w \wed \ast_h({\w^*}) = \la\w,\w\ra,
$ the pairing must be non-degenerate.
\qed

\begin{cor}
If $\Om^\bullet$ has finite dimensional $\del$- and $\adel$-cohomology groups, then 
\bas
H^{(a,b)}_{\del} \simeq (H^{(n-a,n-b)}_{\del})^*, & & H^{(a,b)}_{\adel} \simeq (H^{(n-a,n-b)}_{\adel})^*.
\eas
\end{cor}

\section{Noncommutative K\"ahler Structures}

In this section the definition of a \nc K\"ahler structure is introduced and some of the basic results of classical K\"ahler geometry generalised, most notably the K\"ahler identities. Equality up to scalar multiple of the three Laplacians $\DEL_\exd,$$\DEL_\del$, and $\DEL_{\adel}$,  then follows, implying in turn that  Dolbeault cohomology refines de Rham cohomology. A \nc generalisation of the hard Lefschetz theorem and the $\del \adel$-lemma is then given. The Hermitian structure of $\cpn$ is observed to be K\"ahler, implying that the calculus has cohomology groups of at least classical dimension. Finally, we finish with some spectral calculations and a conjecture about constructing spectral triples for $\bC_q[G/L_S]$.

Throughout this section,  $\Om^{\bullet}$ denotes a covariant $*$-calculus, of total dimension $2n$, over a \qhs $M$. Moreover, $(\Om^{(\bullet,\bullet)},\k)$ denotes  a positive definite covariant Hermitian structure \st the associated integral is closed.

\subsection{K\"ahler Structures and the First Set of K\"ahler Identities}
   
Building on the definition of an Hermitian structure, we define the notion of a  K\"ahler structure. In the classical case this reduces to the fundamental form of a uniquely defined K\"ahler metric \cite[\textsection{3.1}]{HUY}.

\begin{defn}
A {\em K\"ahler structure} for a differential $*$-calculus is an Hermitian structure $(\Om^{(\bullet,\bullet)},\k)$ \st the Hermitian form $\k$ is $\exd$-closed. We call such a $\k$ a  {\em K\"ahler form}.
\end{defn}

Every $2$-form in a $*$-calculus with total dimension $2$ is obviously $\exd$-closed. Hence, just as in the classical case \cite[\textsection{3.1}]{HUY}, \wrt any choice of complex structure, every $\k \in \Om^{(1,1)}$ is a K\"ahler form.

\bigskip

 We now prove the first set of K\"ahler identities. While they follow more or less immediately from the closure of the K\"ahler form, they have important implications throughout the remainder  of the paper.

\begin{lem}\label{KIDSI}
For any K\"ahler structure $(\Om^{(\bullet,\bullet)}, \k)$, we have the following relations
\begin{align*}
[\del,L] = 0, &&  [\adel,L] = 0,       & & [\del^\dagger,\Lam] = 0, &&  [\adel^\dagger, \Lam] = 0. 
\end{align*}
\end{lem}
\demo
By definition a K\"ahler form satisfies $\del \k  = 0$, and so,
\bas
[\del,L](\a) = \del(\k \wed \a) - \k \wed \del \a = \del \k \wed \a + \k \wed \del \a - \k \wed \del \a = 0.
\eas
Analogously, $[\adel, L]  = 0$. The remaining two identities are the  adjoints of the first two.
\qed

\begin{cor}
For every nonzero $\a \in P^k$, there exist unique  forms $\a_0^+, \a_0^- \in P^{k+1}$, $\a_1^+,\a_1^- \in P^{k-1}$ \st
\bal \label{Lefexd}
\del \a = \a_0^+ + L(\a_1^+), & & \adel \a = \a_0^- + L(\a_1^-).
\eal
\end{cor}
\demo
Using the Lefschetz decomposition, $\del \a \in \Om^{k+1}$ can  be written as 
\bas
\del \a =  \sum_{j \geq 0} L^j(\a_j), & & \a_j \in P^{k+1-2j}.
\eas
Since $L$ commutes with $\del$ and $L^{n-k+1} (\a)= 0$, we must have 
$
0 =  \sum_{j \geq 0} L^{n - k + 1+ j}(\a_j). 
$
Moreover, since the Lefschetz decomposition is a direct sum decomposition, we have
\bas
L^{n-k+j+1}(\a_j) = 0, & &  \text{ for all } j \geq 0.
\eas
Now it is only for  $j \leq 1$ that $\a_j$ can be contained in $\ker(L^{n-k+j+1})$. Hence  $\a_j = 0$ for all $j > 2$, and the required identity for $\del$ follows. Uniqueness of $\a_0^+$ is clear. Uniqueness of $\a_1^+$ follows from it being a form of degree at most  $n-1$ and $L$ having trivial kernel in the space of such forms. The proof  for the case of $\adel$ is analogous.
\qed

\subsection{The Second Set of K\"ahler Identities}

In this section we prove the second set of K\"ahler identities and use them to generalise to the \nc setting one of the most important results of K\"ahler geometry, namely that Dolbeault cohomology is a refinement of de Rham cohomology. Throughout this subsection we adopt  the useful convention $L^{j} = 0$, when $j$ is a negative integer.

\begin{lem}\label{LamLCom}
For $\a \in P^k$, it holds that $\Lambda L^j(\a) = [j]_h[n-j-k+1]_h L^{j-1}(\a)$, for $j > 0$.
\end{lem}
\demo
Assuming without loss of generality that $\a \in P^{(a,b)} \sseq P^{k}$, the result follows from 
\bas
~~~~~~~~~~ \Lambda L^j(\a) = & \ast_h\inv L \ast_h L^j(\a) = (-1)^{\frac{k(k+1)}{2}}i^{a-b}\frac{[j]_h!}{[n-j-k]_h!}\ast_h \inv L^{n-j-k+1}(\a) \\
                             = & \, [j]_h [n-j-k+1]_h L^{j-1}(\a). ~~~~~~~~~~~~~~~~~~~~~~~~~~~~~~~~~~~~~~~~~~~~~~~~~~~~ \text{\qed}
\eas

\begin{thm} \label{KIDSII}
The four identities 
\begin{align} \label{KIDS2}
[L,\del^\dagger] = i\ol{\del}, & & [L,\ol{\del}^\dagger] = - i\del, & & [\Lambda,\del] = i\adel^\dagger, & & [\Lambda,\adel] = - i\del^\dagger,
\end{align}
hold in both of the following cases:
\bet
\item the Hodge parameter  is fixed at $h=1$,
\item the domain is restricted to $P^\bullet$ the space of primitive elements.
\eet
\end{thm}
\demo
We begin by finding an explicit description for the action of the left-hand side of the third proposed identity.  For $\a \in P^{(a,b)} \sseq P^k$, and $j \geq 0$, it holds that
\bas
\Lambda \del(L^j(\a)) = & \Lambda (L^j(\del \a)) =  \Lambda L^j(\a^+_0+L(\a^+_1)) =  \Lambda L^j(\a^+_0) + \Lambda L^{j+1}(\a^+_1)\\
= & \, [j]_h[n-j-(k+1)+1]_h L^{j-1}(\a^+_0) \\
   & +  [j+1]_h[n-(j+1)-(k-1)+1]_h L^j(\a^+_1)  \\
= & \, [j]_h[n-j-k]_h L^{j-1}(\a^+_0) +  [j+1]_h[n-j-k+1]_h L^j(\a^+_1). 
\eas
It follows from the above lemma that
\bas
 \del \Lambda(L^j(\a)) = [j]_h[n-j-k+1]_h\big(L^{j-1}(\a^+_0) + L^{j}(\a^+_1)\big).
\eas
Putting these two result together gives
\bas
[\Lambda,\del](L^j(\a))  = & \, \big([j]_h[n-j-k]_h -  [j]_h[n-j-k+1]_h\big)L^{j-1}(\a^+_0)\\
   & +  \big([j+1]_h[n-j-k+1]_h -  [j]_h[n-j-k+1]_h\big)L^j(\a^+_1).
\eas
Moving now to the right-hand side of the proposed identity, we see that
\bas
   i\adel^\dagger(L^j(\a)) 
= & - i \ast_h \del \ast_h(L^j(\a)) \\
= & (-1)^{\frac{k(k+1)}{2}} i^{a-b-1} \frac{[j]_h!}{[n-j-k]_h!} \ast_h \big(L^{n-j-k}(\del(\a))\big)\\
                           = &  (-1)^{\frac{k(k+1)}{2}}i^{a-b-1} \frac{[j]_h!}{[n-j-k]_h!} \ast_h \big(L^{n-j-k}(\a_0^+) + L^{n-j-k+1}(\a^+_1)\big)\\
=   & -  [j]_h L^{j-1}(\a^+_0) + [n-j-k+1]_h L^j(\a^+_1).
\eas
We are now ready to show that the third identity holds in both cases considered above. In the first case, that is when $h = 1$, we have 
\bas
 [\Lambda,\del](L^j(\a)) = & (j(n-j-k) -  j(n-j-k+1))L^{j-1}(\a^+_0)\\
   & + ((j+1)(n-j-k+1) -  j(n-j-k+1))L^j(\a^+_1)\\
 = & - j L^{j-1}(\a_0) + (n-j-k+1)L^j(\a^+_1)\\
 = & \,  i\adel^\dagger(L^j(\a)).
\eas
In the second case, that is when $j=0$, we have 
\bas
 [\Lambda,\del](\a) = [n-k+1]_h (\a_1^+) =  i\adel^\dagger(\a).
\eas

We now move on to the  fourth identity, starting with the case of $h=1$. It follows from the third identity that
\bas
\big(-i\del^*(\w)\big)^* = & \,  i \adel^*(\w^*) = [\Lambda,\del](\w^*) =  \big([\Lambda,\adel](\w)\big)^*, & & \w \in \Om^\bullet.
\eas
Hence,  $[\Lambda,\adel] = - i\del^*$ as required. The second case, that is when $j=0$, is proved analogously using the fact that $P^\bullet$ is closed under the $*$-map.

Finally, we come to the first two identities. For the case of $h=1$, they are  obtained as the adjoints of the first two \wrt the associated inner product. For the case of $j=0$, note that the explicit formulae calculated above for the action of $i \del^\dagger$ and $[\Lambda,\del]$ on $L^j(\a)$ imply that $\del^\dagger,\adel^\dagger, [\Lambda,\del]$, and $[\Lambda,\adel]$, each map $P^\bullet$ to itself.  Since the  Lefschetz decomposition is orthogonal \wrt the inner product, this means that, for $j=0$,  the first and second formulae  can also be obtained by taking adjoints.
\qed


\begin{cor} \label{LaplacianEq}
When  the Hodge parameter is fixed at $h=1$, it holds that
\bas
\del\ol{\del}^\dagger + \ol{\del}^\dagger \del = 0, & &  \del^\dagger \ol{\del} + \ol{\del}\del^\dagger = 0, & & \DEL_{\exd} = 2\DEL_{\del} = 2 \DEL_{\adel}.
\eas
\end{cor}
\demo
The first identity follows from
\bas
i(\del\adel^\dagger + \adel^\dagger\del) & = \del[\Lambda,\del] + [\Lambda,\del]\del = \del\Lambda\del - \del^2\Lambda + \Lambda\del^2 - \del\Lambda\del = 0.
\eas
The second identity is  the operator adjoint of the first.

Moving on to the third identity, we note first that
\begin{align*}
\DEL_\exd = \exd \exd^\dagger + \exd^\dagger \exd & = (\del +  \ol{\del})(\del^\dagger +  \ol{\del}^\dagger) + (\del^\dagger +  \ol{\del}^\dagger)(\del +  \ol{\del})\\
     & = (\del\del^\dagger+\del^\dagger\del) + (\ol{\del} \ol{\del}^\dagger+ \ol{\del}^\dagger \ol{\del})  + (\del \ol{\del}^\dagger +  \ol{\del}^\dagger\del) + (\ol{\del}\del^\dagger + \del^\dagger \ol{\del})\\
     & = \DEL_{\del} + \DEL_{\ol{\del}}.
\end{align*}
It remains to show that $ \DEL_{\del} = \DEL_{\ol{\del}}$:
\begin{align*}
-i \DEL_{\del} =   -i(\del\del^\dagger + \del^\dagger\del) = & \del[\Lambda,\ol{\del}] + [\Lambda,\ol{\del}]\del 
=  \, \del \Lambda \ol{\del} - \del \ol{\del} \Lambda + \Lambda \ol{\del} \del - \ol{\del}\Lambda \del\\
             =  & \, \del\Lambda \ol{\del} - \Lambda \del \adel + \ol{\del} \del \Lambda - \ol{\del} \Lambda \del 
=  [\del,\Lambda]\ol{\del}+\ol{\del}[\del,\Lambda] \\
= & -i\ol{\del}^\dagger\ol{\del} - i\ol{\del}\ol{\del}^\dagger 
=  -i\DEL_{\ol{\del}}. ~~~~~~~~~~~~~~~~~~~~~~~~~~~~~~~~~~~~ \text{\qed}
\end{align*} 


Proportionality of the Laplacians obviously implies equality of harmonic forms:
\bal \label{cohdecomp}
 \H^k_\exd = \bigoplus_{a+b = k} \H^{(a,b)}_{\del} = \bigoplus_{a+b=k} H^{(a,b)}_{\adel}.
\eal  
Hence, Corollary \ref{harmonictoclass} implies the following decomposition of cohomology classes.

\begin{cor} \label{Hdecomp}
De Rham cohomology is {\em refined} by Dolbeault cohomology, which is to say, 
\bal \label{refinement}
H^k_\exd \simeq \bigoplus_{a+b = k} H^{(a,b)}_{\del} \simeq \bigoplus_{a+b = k} H^{(a,b)}_{\adel}.
\eal  
Moreover, the decomposition is independent of the choice of K\"ahler form.
\end{cor}
\demo
We just need to show independence of the decomposition. For $\k'$ another K\"ahler form, denote by $\H^{(a,b)}_{\del}(\k')$ the corresponding space of harmonic forms. For $\w \in \H^{(a,b)}_{\del}$, let $\nu$ be the corresponding element in $\H^{(a,b)}_{\del}(\k')$ \wrt the commutative diagram: 
\begin{align*}
\xymatrix{
  &  \H^{(a,b)}_{\del}  \ar@{^{(}->}[rd] &  \\
H^{(a,b)}_{\del}  \ar[ru]^{\simeq} \ar[rd]_{\simeq} & &  \H^k_\exd \simeq H^k_\exd. \\
 & \H^{(a,b)}_{\del}(\k')  \ar@{_{(}->}[ru]&  \\
}
\end{align*}
We want to show that $\w = \nu + \exd \rho$, for some $\rho \in \Om^{(a-1,b-1)}$. We note first that $\w = \nu + \del \rho'$, for some $\rho \in \Om^{(a-1,b)}$. Moreover, $\del \rho'$ is $\exd$-closed because $\exd(\del \rho') = \exd(\w - \nu) = 0$. By Hodge decomposition \wrt $\exd$, this means $\del \rho'$ is the sum of a harmonic form and a  $\exd$-exact form. But  Corollary \ref{harmonictoclass}  tells us that $\del \rho'$ is contained in a complementary subspace to $\H^\bullet$. Hence, it must be $\exd$-exact and independence of the decomposition follows.
\qed

\subsection{Harmonic Forms and the Hodge Parameter}

We begin by showing that the Lefschetz and dual Lefschetz operators commute with the Laplacian $\DEL_\exd$, and hence that they induce operators on the space of harmonic forms. 

\begin{lem}\label{Lcommutes}
When the Hodge parameter is fixed at $h=1$, 
\bas
[L,\DEL_\exd] = [\Lambda,\DEL_\exd] = 0.
\eas
\end{lem}
\demo
Using the K\"ahler identity $L\adel^\dagger = \adel^\dagger L -i\del$, and proportionality of the Laplacians, we see that
\bas
\frac{1}{2} L \DEL_{\exd}  = &   L \DEL_{\adel} =  L(\adel\adel^\dagger + \adel^\dagger \adel) = \adel L \adel^\dagger + (\adel^\dagger L -i\del)\adel \\
= &\,  \adel (\adel^\dagger L -i\del) + \adel^\dagger \adel L + i\adel \del =  (\adel\adel^\dagger + \adel^\dagger \adel) L \\
= &  \DEL_{\del} L =  \frac{1}{2} \DEL_{\exd} L.
\eas
The second relation is the adjoint of the first. 
\qed

Up to this point we have avoided the question of whether the the space of harmonic forms depends on the Hodge parameter. The following lemma and its corollary provides an answer to this question.

\begin{lem} \label{HarmSumL}
When the Hodge parameter is fixed at $h=1$, a form $\w$  with Lefschetz decomposition $\w = \sum_{j \geq 0} L^j(\a_j)$, for $\a_j \in P^{\bullet}$, is harmonic \iff $\a_j$ is $\exd$-closed, for all $j$.
\end{lem}
\demo 
Lemma \ref{LamLCom} implies that
\bas
 \Lambda^m(\w) = &  \Lambda^m \bigg( \sum_{j=1}^m L^j(\a_j)\bigg) =  \Lambda^m L^m(\a_m)   =  \bigg(\prod_{j=1}^m j(n-j-k+1) \bigg) \a_{m}.
\eas
By the above lemma, $\H^\bullet$ is closed under $L$ and $\Lambda$, and so, if $\w \in \H^\bullet$ then $\a_m \in \H^\bullet$. This in turn implies that
\bas 
\w - L^m(\a_{m}) = \sum_{j=1}^{m-1} L^j(\a_{j}) \in \H^\bullet.
\eas
Repeated applications of this argument show that if $\w \in \H^\bullet$  then $\a_j \in \H^\bullet$, for all $j$. The converse follows from the closure of $\H^\bullet$ under $L$ and $\Lambda$, and so, $\w$ is harmonic \iff $\a_j \in \H^\bullet$, for all $j$.   

It remains to show that a primitive form is harmonic \iff it is $\exd$-closed. Clearly, we need only show that $\exd$-closure implies harmonicity. This follows from  Lemma \ref{harmcap} and the fact that, for $\a \in P^{(a,b)} \sseq P^k$, we have
\bas
~~~~~~~~~~~~~~~~~~~ \exd^*\a = & \ast_1 \exd \ast_1(\a) = (-1)^{\frac{k(k+1)}{2}} i^{a-b} \frac{1}{(n-k)!} \ast_1 \exd \big(L^{n-k}(\a)\big) \\
=  & \, (-1)^{\frac{k(k+1)}{2}} i^{a-b} \frac{1}{(n-k)!} \ast_1  L^{n-k}(\exd \a) = 0. ~~~~~~~~~~~~~~~~~~~~~~~~~~~~~~~~ \text{\qed}
\eas

The proposition shows us that, for $h=1$, the space of harmonic forms is completely determined by the $\exd$-closed primitive forms.  The following corollary tells us that this is also the case when $h \neq 1$.

\begin{cor}
For any choice of Hodge parameter $h \in \bR_{>0}$, it holds that
\bas
\H^\bullet_{\exd} = \H^\bullet_{\del} = \H^\bullet_{\adel} = \spn_{\bC}\big\{L^j(\a) \,|\, j \in\bN_0, \a \in P^\bullet \cap \ker(\exd)\big\}.
\eas
\end{cor}
\demo
The fact that  $L^j(\a) \in \H^\bullet_{\exd}, \H^\bullet_{\del}$, and $\H^\bullet_{\adel}$, for any value of $h$, is shown just as in the $h=1$ case. Hence 
\bas
\spn_{\bC}\big\{L^j(\a) \,|\, j \in\bN_0, \a \in P^\bullet \cap \ker(\exd)\big\} \sseq \H^\bullet_{\exd}, \H^\bullet_{\del}, \H^\bullet_{\adel}.
\eas  
Since the cohomology groups $H^\bullet_{\exd}, H^\bullet_{\del}$, and $H^\bullet_{\adel}$, are defined independently of $h$, Corollary \ref{harmonictoclass} now implies that these inclusions are equalities.  
\qed

\subsection{The Hard Lefschetz Theorem}

Lemma  \ref{Lcommutes}  and Corollary \ref{lapcommiso} imply that $L$ and $\Lambda$ induce maps on $H^\bullet$. This  allows us to make the following definition.

\begin{defn}
For a K\"ahler structure, the {\em $(a,b)$-primitive cohomology group} is the vector space
\bas
H^{(a,b)}_{\text{prim}} := \ker\Big(L^{n-(a+b)+1}:H^{(a,b)} \to H^{(n-b+1,n-a+1)}\Big).
\eas
Moreover, we denote $H^{k}_{\text{prim}} := \bigoplus_{a+b=k} H^{(a,b)}_{\text{prim}} $.
\end{defn}

This definition, together with Proposition \ref{HarmSumL} and Lemma \ref{harmonictoclass}, gives us the following \nc generalisation of the classical hard Lefschetz theorem \cite[Proposition 3.3.13]{HUY}. As a corollary we prove a generalisation Corollary \ref{h2ncor} to the case of $H^{2k}$, for all $k = 0, \ldots, n$.

\begin{thm}
Let $(\Om^{(\bullet,\bullet)},\exd)$ be a K\"ahler structure, then it holds that
\bet
\item $L^{k}: H^{n-k} \to H^{n+k}$ is an isomorphism, for $k = 0, \ldots, n$,
\item $H^k \simeq \bigoplus_{i\geq 0} L^i H^{(a,b)}_{\text{prim}}$.
\eet
\end{thm}

\begin{cor}\label{cohombdd}
For a covariant differential $*$-calculus endowed with a  covariant  K\"ahler structure,  it holds that
\bas
\dim(H^{2k}) = \dim(H^0) \geq 1, & & \text{ for all } k = 0, \ldots, n.
\eas
\end{cor}
\demo
This is a direct consequence of the second statement of above theorem and the fact that $\exd 1 = 0$.
\qed

\subsection{The $\del \adel$-Lemma}

We finish our study of the general theory of K\"ahler structures with a result known in the classical case as the $\del \adel$-lemma.  While it may look like an innocent technical result, in the classical case it is crucial for many important results, such as formality for K\"ahler manifolds \cite[\textsection 3.A]{HUY}.

\begin{lem} Let $\Om^\bullet$ be a covariant differential  $*$-calculus admitting a covariant K\"ahler structure. Then for a $\exd$-closed form $\w \in \Om^{(a,b)}$, the following conditions are equivalent:
\begin{enumerate}
\item $\w$ is  $\exd$-exact,
\item $\w$ is $\del$-exact,
\item $\w$ is  $\adel$-exact,
\item $\w$ is  $\del \adel$-exact.
\end{enumerate}
\end{lem}
\demo
We will prove the theorem by introducing a fifth equivalent condition: $\w$ is  orthogonal to $\H^{(a,b)}$ for some choice of K\"ahler form. 

Using Hodge decomposition, we see the fifth condition is implied by any of the other four  conditions. Moreover, the fourth condition implies both the first, second, and third conditions. Thus, it suffices to show the fifth condition implies the fourth.

Since by assumption $\w$ is $\exd$-closed (and hence $\del$-closed) and orthogonal to the space of harmonic forms, then Hodge decomposition with respect to $\del$ yields $\w = \del \nu$, for some $\nu \in \Om^{(a-1,b)}$. Applying Hodge decomposition \wrt $\adel$ to  $\nu$ yields $\nu = \adel \nu' + \adel^\dagger \nu'' + \nu'''$, for some harmonic $\nu'''$. Returning to original form $\w$, we now see that $\w = \del \adel \nu' + \del \adel^\dagger \nu''$. By assumption $\adel \w = 0$, and so, fixing $h=1$, Corollary \ref{LaplacianEq}  implies that
\bas
0 = \adel \w = \adel \del \adel \nu' + \adel \del \adel^\dagger \nu' = - \adel \adel^\dagger \del \nu'.
\eas
Since $0 = \big< \adel\adel^\dagger \del \nu',\del \nu'\big> = \big<\adel^\dagger \del \nu', \adel^\dagger \del \nu'\big>$, this means that $\adel^\dagger \del \nu' = 0$. Thus $\w = \del \adel \nu'$.
\qed

\subsection{The \hk Calculus}

The next result follows directly from Lemma 4.17 and Proposition 4.19.

\begin{lem} The Hermitian structure $(\Om^{(\bullet,\bullet)},\k)$ for $\ccpn$ is a K\"ahler structure.
\end{lem}

The operator  $\adel +\adel^\dagger$ is a direct $q$-deformation of the  Dirac--Dolbeault operator of $\cpn$. Deformations of this operator have previously appeared in the literature \cite{DDCPN} in the context of spectral triples \cite[Chapter 10]{VAR}.  As an initial investigation of the spectrum of  \linebreak $\adel +\adel^\dagger$, we calculate the first non-zero eigenvalue of the Laplacian $\DEL_{\adel}$,  for $\bC_q[\bC P^1]$ in the following lemma and corollary.

\begin{lem}
For  $X,Y:\bC_q[\bC P^1]  \to \bC$  the linear functionals uniquely defined by 
$
X(m)e^+ + Y(m)e^- := [m^+],
$
it holds that 
\bas
\DEL_{\adel}(m) = - X(m\2)Y(m\3) m\1, & & m \in \bC_q[\bC P^1] = \Om^{(0,0)}.
\eas
\end{lem}
\demo
Suppressing explicit reference to $\unit$, we see that 
\bas
\DEL_{\adel}(m) = & \ast_q \circ \, \del \circ \ast_q \circ \adel(m) \\
= &  \ast_q \circ \, \del \circ \ast_q\big(m\1 X(m\2) \oby e^-_1 \big)\\
= &  -i \ast_q \circ \, \del \big(m\1 X(m\2) \oby e^-_1 \big)\\
= &  -i \ast_q  \big(m\1 Y(m\2)X(m\3) \oby e^+_1 \wed e^-_1 \big)\\
= &  - m\1 Y(m\2)X(m\3),
\eas
where we have used \cite[Lemma 5.5]{MMF2} to calculate the actions of $\del$ and $\adel$.
\qed

\begin{cor}
It holds that, 
$
\DEL_{\adel}(z_{ij}) =  q [2]_{q} z_{ij}$, for $i \neq j.
$
\end{cor}
\demo
From the above lemma, we have that 
\bas
\DEL_{\adel}(z_{ij})  = & - X((z_{ij})\2)Y((z_{ij})\3) (z_{ij})\1\\
  = &  \, - \sum_{a,b, x,y=1}^2 X(u^a_bS(u^y_x))Y(u^b_1S(u^1_y))u^i_a S(u^x_j).
\eas
Using the  formulae presented in  \cite[Proposition 3.3]{MMF2}, it is easily calculated that the scalars $X(u^a_bS(u^y_x))Y(u^b_1S(u^1_y))$ are non-zero only in the following cases
\bas
X(u^1_1S(u^2_1))Y(u^1_1S(u^1_2)) = - q^2 , & & X(u^2_1S(u^2_2))Y(u^1_1S(u^1_2)) = 1.
\eas
Hence  $\DEL_{\adel}(z_{ij}) =  q^2 u^i_1S(u^1_j) -  u^i_2S(u^2_j)$. Finally, our assumption  that $i \neq j$ implies that $u^i_2S(u^2_j) = - u^i_1S(u^1_j)$, and so, 
\bas
~~~~~~~~~~~~~~~~~~~~  \DEL_{\adel}(z_{ij}) = & \, q^2 u^i_1S(u^1_j) +   u^i_1S(u^1_j)  = q (q + q\inv)z_{ij}  
=  q [2]_q z_{ij}.  ~~~~~~~~~~~~~~~~ \text{\qed}
\eas

Recall now $(\Om^{(\bullet,\bullet)}, \k)$ the conjectured Hermitian structure  for $\bC_q[G/L_S]$ introduced in \textsection 4.5. Using a direct generalisation of Lemma 4.17, it can be shown that $\k$ must be $\exd$-closed, and so, we have the following lemma. 

\begin{lem}
If the pair $(\Om^{(\bullet,\bullet)}, \k)$ is an Hermitian structure  for $\bC_q[G/L_S]$, then it is a K\"ahler structure.
\end{lem}

We finish with a conjecture about the completion of these conjectured K\"ahler structures to spectral triples for the irreducible quantum flag manifolds. The case of $\cpn$ is treated in \cite{DOS}.

\begin{conj}
Denoting by  $L^2(\Om^{(0,\bullet)})$ the completion of the subcomplex $\Om^{(0,\bullet)}$ of the Heckenberger--Kolb calculus of $\bC_q[G/L_S]$  \wrt the inner product associated to $\k$, a spectral triple is given by 
\bas
(\bC_q[G/L_S], L^2(\Om^{(0,\bullet)}), \adel + \adel^\dagger).
\eas
\end{conj}


\bigskip

Instytut Matematyczny, Polskiej Akademii Nauk, ul. \'Sniadeckich 8, 00-656 Warszawa, Poland

{\em e-mail}: \tt{robuachalla@impan.pl}

\end{document}